\documentclass[12pt,reqno]{amsart}
\usepackage{amssymb, amsfonts, amsbsy, latexsym, epsfig, color}

\newtheorem{Thm}{Theorem}[section]

\newtheorem{lemma}[Thm]{Lemma}

\newtheorem{proposition}[Thm]{Proposition}
\newtheorem{definition}[Thm]{Definition}

\newtheorem{example}[Thm]{Example}
\newtheorem{theorem}[Thm]{Theorem}

\newcommand{\bitem}{\begin{itemize}}
\newcommand{\eitem}{\end{itemize}}
\newcommand{\benum}{\begin{enumerate}}
\newcommand{\eenum}{\end{enumerate}}
\newcommand{\beq}{\begin{equation}}
\newcommand{\eeq}{\end{equation}}

\newcommand{\norm}[1]{\|#1\|}

  \newcommand{\R}{\mathbb{R}}

 \newcommand{\Z}{\mathbb{Z}}
 
 \DeclareMathOperator*{\esssup}{ess\,sup}
\DeclareMathOperator*{\essinf}{ess\,inf}
\DeclareMathOperator*{\supp}{supp}

\def\ZZ{\mathbb{Z}}
\def\RR{\mathbb{R}}
\def\ZZ{\mathbb{Z}}
\def\CC{\mathbb{C}}

\newcommand{\bZ}{{\mathbb Z}}
\newcommand{\bR}{{\mathbb R}}

\newcommand{\bS}{{\mathbb S}}

\def\cC{{\mathcal{C}}}

\def\cE{{\mathcal{E}}}

\def\cR{{\mathcal{R}}}

\def\cSH{{\mathcal{S}\mathcal{H}}}

\newcommand{\wq}[1]{{\color{black}{#1}}}

\begin{document}

\title{Construction of Compactly Supported Shearlet Frames}

\author[P. Kittipoom]{Pisamai Kittipoom}
\address{Faculty of Science, Prince of Songkla University, Hat Yai, Songkhla 90112, Thailand}
\email{pisamai.k@psu.ac.th}

\author[G. Kutyniok]{Gitta Kutyniok}
\address{Institute of Mathematics, University of Osnabr\"uck, 49069 Osnabr\"uck, Germany}
\email{kutyniok@math.uni-osnabrueck.de}

\author[W.-Q Lim]{Wang-Q Lim}
\address{Institute of Mathematics, University of Osnabr\"uck, 49069 Osnabr\"uck, Germany}
\email{wlim@math.uni-osnabrueck.de}

\thanks{P.K. acknowledges support from the Faculty of Science, Prince of Songkla University. G.K. and W.-Q L.
would like to thank Wolfgang Dahmen, Chunyan Huang, Demetrio Labate, Christoph Schwab,
and Gerrit Welper for various discussions on related topics. G.K. and W.-Q L. acknowledge
support from DFG Grant SPP-1324, KU 1446/13-1.}

\begin{abstract}
Shearlet tight frames have been extensively studied during the last years due to their optimal
approximation properties of cartoon-like images and their unified treatment of the continuum
and digital setting. However, these studies only concerned shearlet tight frames generated by
a band-limited shearlet, whereas for practical purposes compact support in spatial domain is
crucial.

In this paper, we focus on cone-adapted shearlet systems which
-- accounting for stability questions -- are associated with a general irregular set of parameters.
We first derive sufficient conditions for such cone-adapted irregular shearlet systems to form a
frame and provide explicit estimates for their frame bounds. Secondly, exploring these results
and using specifically designed wavelet scaling functions and filters,
we construct a family of cone-adapted shearlet frames consisting of compactly supported
shearlets. For this family, we derive
estimates for the ratio of their frame bounds and prove that they provide optimally
sparse approximations of cartoon-like images.
\end{abstract}

\keywords{Curvilinear discontinuities, edges, filters, nonlinear approximation,
optimal sparsity, scaling functions, shearlets, wavelets}

\subjclass{Primary: 42C40; Secondary: 42C15, 65T60, 65T99, 94A08}

\maketitle

\section{Introduction}

Over the last 20 years wavelets have established themselves as a key methodology for efficiently representing
signals or operators with applications ranging from more theoretical tasks such as adaptive schemes for solving elliptic
partial differential equations to more practically tasks such as data compression. It is fair to say that nowadays
wavelet theory can be regarded as an essential area in applied mathematics. A major breakthrough in this field was
achieved through Daubechies in 1988 by introducing compactly supported wavelet orthonormal bases with favorable vanishing
moment and support properties \cite{Dau88}.

Recently, it was proven by Cand\'{e}s and Donoho in \cite{CD04} that wavelets do not perform optimally
when representing and analyzing anisotropic features in multivariate data such as singularities concentrated
on lower dimensional embedded manifolds, where their focus was on the 2-dimensional situation. This observation
initiated a flurry of activity to design a representation system which provides optimally sparse approximations
of so-called cartoon-like images while having as favorable theoretical and computational properties as
wavelet systems (see, e.g., \cite{CD04,DV05}). In 2005, {\em shearlets} were developed by Labate,
Weiss, and two of the authors in \cite{LLKW05} (see also \cite{GKL06}) as the first directional representation
system with allows a unified treatment of the continuum and digital world similar to wavelets, while
providing (almost) optimally sparse approximations within a cartoon-like model \cite{GL07} -- `almost'
in the sense of an additional log-factor which is customarily regarded as negligible.

Several constructions of discrete shearlet frames are already known to date, see \cite{GKL06,KL07,DKST09,Lim09,KKL10}.
However, taking applications into account, spatial localization of the analyzing elements of the encoding system is of
uttermost importance both for a precise detection of geometric features as well as for a fast decomposition
algorithm. But, one must admit that non of the previous approaches encompasses this crucial case. Hence the
main goal of this paper is to provide a comprehensive analysis of cone-adapted -- the variant of shearlet systems
feasible for applications -- discrete shearlet frames encompassing in particular compactly supported shearlet generators.
Our contribution is two-fold: We firstly provide sufficient conditions for a cone-adapted irregular shearlet
system to form a frame for $L^2(\RR^2)$ with explicit estimates for the ratio of the associated frame bounds;
and secondly, based on these results, we introduce a class of cone-adapted compactly supported shearlet frames,
which are even shown to provide (almost) optimally sparse approximations of cartoon-like images, alongside with
estimates for the ratio of the associated frame bounds.


\subsection{Shearlet Systems}

Referring to the detailed introduction of shearlet systems in Section \ref{sec:shearletsystems}, we
allow us here to just briefly mention those main ideas and results, which are crucial for an intuitive understanding.

Shearlet systems are designed to efficiently encode anisotropic features. In order to achieve optimal sparsity,
shearlets are scaled according to a parabolic scaling law, thereby exhibiting a spatial footprint of size $2^{-j}$
times $2^{-j/2}$. They parameterize directions by slope encoded in a shear matrix, since
choosing shear rather than rotation is in fact decisive for a unified treatment of the continuum and
digital setting.  We refer for more details to \cite{GKL06,KL09} for the continuum theory and
\cite{DKS08,ELL08,Lim09} for the digital theory.

Shearlet systems come in two ways: One class being generated by a unitary representation of the shearlet
group and equipped with a particularly `nice' mathematical structure, however causes a biasedness
towards one axis, which makes it unattractive for applications; the
other class being generated by a quite similar procedure however restricted to a horizontal and vertical
cone in frequency domain, thereby ensuring an equal treatment of all directions.
For both cases, the continuous shearlet systems are associated with a 4-dimensional parameter space
consisting of a scale parameter measuring the resolution, a shear parameter measuring the orientation,
and a translation parameter measuring the position of the shearlet. A sampling of this parameter
space leads to discrete shearlet systems, and it is obvious that the possibilities for this are numerous.
A canonical sampling approach -- using dyadic sampling -- leads to so-called regular
shearlet systems. However, lately, questions from sampling theory and stability issues brought the study
of irregular shearlet systems to the researcher's attention.


\subsection{Necessity of Compactly Supported Shearlets}

Shearlet theory has already impacted various applications for which the sparse
encoding or analysis of anisotropic features is crucial, and we refer to results on denoising \cite{Lim09},
edge detection \cite{GLL09}, or geometric separation \cite{DK09,DK10}. However, similar as in
wavelet analysis, a significant improvement of the applicability of shearlets is possible,
if compactly supported shearlet frames would be made available. Let us discuss two
examples to illustrate the necessity of such a study.

\textit{Imaging Science.} In computer vision, edges were detected as those features governing an image
while separating smooth regions in between. Thus imaging tasks typically concern an especially careful
handling of edges, for instance, avoiding to smoothen them during a denoising process. Therefore,
when exploiting a decomposition in terms of a representation system for such tasks, a superior
localization in spatial domain of the analyzing elements is in need. This would then allow a very
precise focus on edges, thereby reducing or even avoiding various artifacts in the decomposition
and analysis process. Secondly, due to the increased necessity to process very large data sets,
fast decomposition algorithms are crucial. Learning from the experience with structured transforms
such as the fast wavelet transform, it is conceivable that a representation system consisting of
compactly supported elements will provide a significant gain in complexity.

\textit{Partial Differential Equations.} Adaptive schemes for solving elliptic partial differential
equations using wavelet decompositions have turned out to be highly beneficial both in theory and applications
(cf. \cite{Dah97}). However, hyperbolic partial differential equations exhibit shock fronts, which
are not optimally sparsely encoded by wavelets due to their isotropic footprints. For a directional
representation system such as shearlets to be feasible as an encoding method for an adaptive scheme,
one main obstacle is the sparsity of the associated stiffness matrix. However, with good localization
properties such as compact support, this obstacle might be overcome. Another problem is the handling
of boundary conditions on a bounded domain, which seems also attackable by having compactly supported
shearlet frames at hand.
Last, but not least, again a fast decomposition algorithm is of importance, and we refer to the
previously described application for a reasoning why having compact support will make a tremendous difference here.


\subsection{Previous Work on Shearlet Constructions}

Let us now discuss the history and the state-of-the-art of shearlet constructions. The
first class of shearlets, which were shown to generate a tight frame, were band-limited
with a wedge-like support in frequency domain specifically adapted to the shearing operation (see \cite{LLKW05,GKL06}).
This particular class of cone-adapted shearlet frames was already extensively explored, for
instance, for analyzing sparse approximation properties of the associated shearlet frames
\cite{GL07}.

Shortly afterwards, a different avenue was undertaken in \cite{KL07}, where a first attempt
was made to derive sufficient conditions for the existence of irregular shearlet frames, which,
however, could not decouple from the previous focus on band-limited shearlet generators. In
addition, this result was stated for shearlet frames which arose directly from a group representation.
In some sense, this path was continued in \cite{DKST09}, where again sufficient conditions for
this class of irregular shearlet frames were studied. Let us also mention at this point that a
more extensive study of these systems was performed in \cite{KKL10} with a focus on necessary
conditions and a geometric analysis of the irregular parameter set.

Now returning to the situation of cone-adapted
shearlet frames, a particular interesting study was recently done in \cite{Lim09}, where a fast
decomposition algorithm for {\em separable} shearlet frames -- which are non-tight -- was
introduced, however without analyzing frame or sparsity properties. We should mention that
it was this work which gave us the intuition of taking a viewpoint of separability for the
construction of compactly supported shearlet frames in this paper. The question
of optimal sparsity of a large class of cone-adapted compactly supported shearlet frames
was very recently solved by two of the authors in \cite{KL10}. However, no sufficient conditions
for or construction of cone-adapted compactly supported shearlet frames with estimates for the
frame bounds nor with sparsity analyses are known to date. This is an open question even for
other directional representation systems such as the perhaps previously most well-known curvelet
system \cite{CD04}.


\subsection{Our Contribution}

Our contribution in this paper is two-fold. Firstly, Theorem \ref{theo:general} provides sufficient
conditions on the irregular set of parameters and the generating shearlet to form a  cone-adapted
shearlet frame alongside with estimates for the ratio of the frame bounds. This result greatly
extends the result in \cite{KL07} by, in particular, encompassing compactly supported shearlet
generators. Secondly, exploring these results and using specifically designed wavelet scaling
functions and filters, we explicitly construct a family of cone-adapted shearlet frames consisting
of compactly supported shearlets with estimates for their frame bounds in two steps: For functions
supported on the horizontal cone in frequency domain and then for functions in $L^2(\RR^2)$ (see
Theorems \ref{theo:compact} and \ref{theo:completeframe}). This construction is the first of its
kind in shearlet theory, and it is worth mentioning that, to the author's knowledge, even in wavelet theory there does not
exist a comparable one for compactly supported wavelets. Theorem \ref{theo:sparse} then proves
that this family does even provide (almost) optimally sparse approximations of cartoon-like images.


\subsection{Outline}

In Section \ref{sec:shearletsystems}, we introduce the necessary definitions and notation for
cone-adapted irregular shearlet systems. In particular, in Subsections \ref{subsec:coneshearlets}
and \ref{subsec:irregularsampling} we present the viewpoint of deriving cone-adapted irregular
shearlet systems from sampling the parameters of cone-adapted continuous shearlet systems for
the first time. Subsection \ref{subsec:class} is then concerned with introducing a large class of
cone-adapted irregular shearlet systems, coined feasible shearlet systems, which will be the
focus of this paper. Sufficient conditions for such feasible shearlet systems to form a frame
and explicit estimates for their frame bounds are then discussed in Section \ref{sec:generalsufficient}.
In Section \ref{sec:compactsupport}, these results are explored to construct a family of
cone-adapted regular compactly supported shearlet frames with estimates for the ratio of the
associated frame bounds (Subsections \ref{subsec:compact} and \ref{subsec:generalRR}).
In Subsection \ref{subsec:sparse}, it is finally proven that
this family does even provide (almost) optimally sparse approximations of cartoon-like images. Some very
technical, lengthy proofs are deferred to Section \ref{sec:proofs}.


\section{Regular and Irregular Shearlet Systems}
\label{sec:shearletsystems}

\subsection{Cone-Adapted Shearlets}
\label{subsec:coneshearlets}

Shearlets are highly anisotropic representation systems, which optimally sparsify $C^2(\bR^2)$-functions
apart from $C^2$-discontinuity curves. Their main advantage over other proposed directional
representation systems is the fact that they provide a unified treatment of the continuum
and digital world. Today, shearlets are utilized for various applications, and we refer to \cite{GLL09,DK09,DK10}
as some examples.

Introduced in \cite{GKL06}, discrete shearlet systems -- shearlet systems for
$L^2(\bR^2)$ with discrete parameters -- exist in two different variants: One coming directly from a
group representation of a particular semi-direct product, the so-called shearlet group, and the other being adapted to a cone-like
partitioning of the frequency domain. In fact, the second type of shearlet systems are the ones
exhibiting the favorable property of treating the continuum and digital setting uniformly similar to
wavelets, and are the ones relevant for applications. Hence those are the discrete shearlet systems we focus on in this paper.

We now introduce cone-adapted discrete regular shearlet systems -- previously also referred to as `shearlets
on the cone' -- as a special sampling of the (cone-adapted) continuous shearlet systems, which were
exploited in \cite{KL09}. Since all previous papers consider discrete regular shearlet systems as
sampling the parameters of continuous shearlet systems refer to the first type of shearlet systems
arising from a group representation (see, e.g., \cite{KL07,DKST09}), this is the first time this view
is taken and hence is presented in a more elaborate manner.


\subsubsection{(Cone-Adapted) Continuous Shearlet Systems}


Shearlets are scaled according to a parabolic scaling law encoded in the {\em parabolic scaling matrices}
$A_a$ or $\tilde{A}_{a}$, $a > 0$, and exhibit directionality by parameterizing slope encoded
in the {\em shear matrices} $S_s$, $s \in \RR$, defined by
\[
A_{a} =
\begin{pmatrix}
  a & 0\\ 0 & \sqrt{a}
\end{pmatrix}
\qquad\mbox{or}\qquad
\tilde{A}_{a} =
\begin{pmatrix}
  \sqrt{a} & 0\\ 0 & a
\end{pmatrix}
\]
and
\[
S_s = \begin{pmatrix}
  1 & s\\ 0 & 1
\end{pmatrix},
\]
respectively.

\begin{figure}[ht]
\begin{center}
\includegraphics[height=1.4in]{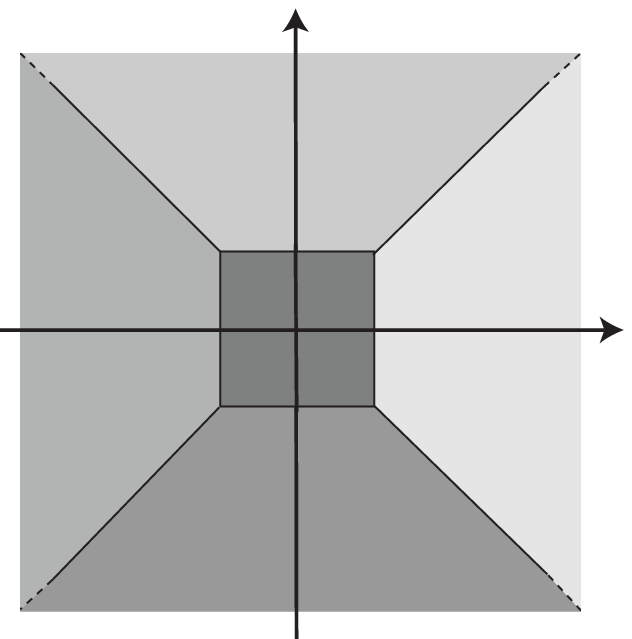}
\put(-33,58){\footnotesize{$\cC=\cC_1$}}
\put(-70,80){\footnotesize{$\cC_2$}}
\put(-88,30){\footnotesize{$\cC_3$}}				
\put(-50,52){\footnotesize{$\cR$}}
\put(-45,15){\footnotesize{$\cC_4$}}
\put(-60,-17){(a)}
\hspace*{4cm}
\includegraphics[height=1.4in]{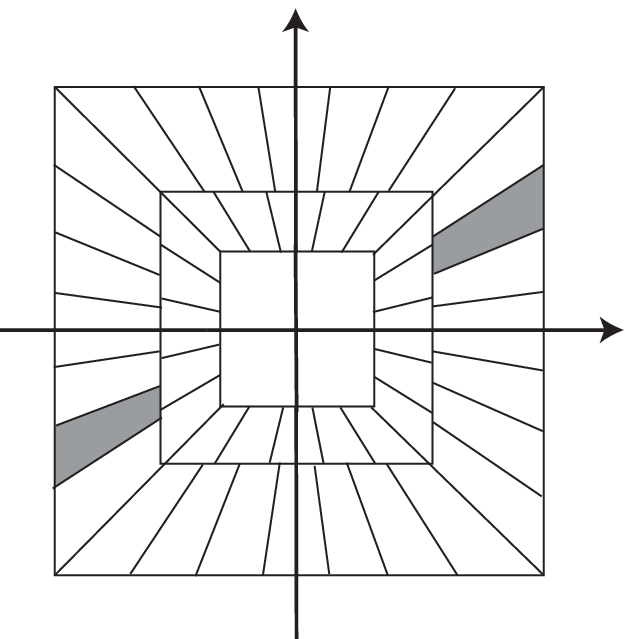}
\put(-60,-17){(b)}
\end{center}
\caption{(a) The cones $\cC_1$ -- $\cC_4$ and the centered rectangle $\cR$ in frequency domain.
(b) The tiling of the frequency domain induced by discrete shearlets with the support of one
shearlet from Example \ref{exa:classicalshearlet} exemplarily highlighted.}
\label{fig:tilings}
\end{figure}

To ensure an (almost) equal treatment of the different slopes, which is
evidently of uttermost importance for practical applications, we partition the frequency
plane into the following four cones $\cC_1$ -- $\cC_4$:
\begin{equation*}
\cC_\iota = \left\{ \begin{array}{rcl}
\{(\xi_1,\xi_2) \in \bR^2 : \xi_1 \ge 1,\, |\xi_2/\xi_1| \le 1\} & : & \iota = 1,\\
\{(\xi_1,\xi_2) \in \bR^2 : \xi_2 \ge 1,\, |\xi_1/\xi_2| \le 1 \} & : & \iota = 2,\\
\{(\xi_1,\xi_2) \in \bR^2 : \xi_1 \le -1,\, |\xi_2/\xi_1| \le 1\} & : & \iota = 3,\\
\{(\xi_1,\xi_2) \in \bR^2 : \xi_2 \le -1,\, |\xi_1/\xi_2| \le 1\} & : & \iota = 4,
\end{array}
\right.
\end{equation*}
and a centered rectangle
\[
\cR = \{(\xi_1,\xi_2) \in \bR^2 : \|(\xi_1,\xi_2)\|_\infty < 1\}.
\]
For an illustration, we refer to Figure \ref{fig:tilings}(a).

The rectangle $\cR$ corresponds
to the low frequency content of a signal, which is customarily represented by translations
of some scaling function. Anisotropy now comes into play when encoding the high frequency content
of a signal, which corresponds to the cones $\cC_1$ -- $\cC_4$, where the cones $\cC_1$ and $\cC_3$
as well as $\cC_2$ and $\cC_4$ are treated separately as can be seen in the following definition.

\begin{definition}
\label{defi:contshearlets}
The {\em (cone-adapted) continuous shearlet system} $\cSH(\phi,\psi,\tilde{\psi})$ generated by a {\em scaling function}
$\phi \in L^2(\mathbb{R}^2)$ and {\em shearlets} $\psi, \tilde{\psi} \in L^2(\mathbb{R}^2)$ is defined by
\[
\cSH(\phi,\psi,\tilde{\psi}) = \Phi(\phi) \cup \Psi(\psi) \cup \tilde{\Psi}(\tilde{\psi}),
\]
where
\[
\Phi(\phi) = \{\phi_t  = \phi(\cdot-t) : t \in \RR^2\},
\]
\[
\Psi(\psi) = \{\psi_{a,s,t}=  a^{-\frac34} \psi(A_a^{-1}S_s^{-1}(\, \cdot \, -t))
: a \in (0,1],\, s \in [-(1+a^{1/2}),1+a^{1/2}],\, t \in \mathbb{R}^2\},
\]
and
\[
\tilde{\Psi}(\tilde{\psi}) = \{\tilde{\psi}_{a,s,t}=  a^{-\frac34} \tilde{\psi}(\tilde{A}_a^{-1}S_s^{-T}(\, \cdot \, -t))
: a \in (0,1],\, s \in [-(1+a^{1/2}),1+a^{1/2}],\, t \in \mathbb{R}^2\},
\]
where $a$ indexes scale, $s$ indexes shear, and $t$ indexes translation. Setting
\[
\bS_{cone}=\{(a,s,t) : a \in (0,1],\, s \in [-(1+a^{1/2}),1+a^{1/2}],\, t \in \mathbb{R}^2\},
\]
the associated {\em (Cone-Adapted) Continuous Shearlet Transform} $\mathcal{S}\mathcal{H}_{\phi,\psi,\tilde{\psi}}
f : \RR^2 \times \bS_{cone}^2 \to \mathbb{C}^3$ of some function $f \in L^2(\RR^2)$ is given by
\[
\mathcal{S}\mathcal{H}_{\phi,\psi,\tilde{\psi}} f(t',(a,s,t),(\tilde{a},\tilde{s},\tilde{t}))
= (\langle f,\phi_{t'}\rangle,\langle f,\psi_{a,s,t}\rangle,\langle f,\tilde{\psi}_{\tilde{a},\tilde{s},\tilde{t}}\rangle).
\]
\end{definition}

Notice that in \cite{KL09} only special generators $\psi$ were considered, whereas here we state
the definition for all $\psi \in L^2(\bR^2)$; see also the interesting extension \cite{Gro09b}.


\subsubsection{(Cone-Adapted) Discrete Shearlet Systems}

A discretization of the shearlet parameters is typically achieved through sampling of the parameter set of
(cone-adapted) continuous shearlet systems, where $\RR^2$ is sampled as $c_1\ZZ^2$ and, for the horizontal cone
$\cC_1 \cup \cC_3$, $\bS_{cone}$ is discretized as
\beq \label{eq:regularsampling}
\{(2^{-j},k2^{-j/2},S_{k2^{-j/2}}A_{2^{-j}}M_cm : j \ge 0,\, k \in
\{-\lceil{2^{j/2}}\rceil,\ldots,\lceil{2^{j/2}}\rceil\},\, m \in \bZ^2\},
\eeq
where, for $c=(c_1,c_2) \in (\RR^+)^2$, $M_c$ denotes the sampling matrix
\[
M_c = \begin{pmatrix}
  c_1 & 0\\ 0 & c_2
\end{pmatrix}.
\]
For the vertical cone $\cC_2 \cup \cC_4$, $\bS_{cone}$ is discretized in a similar way, now using the sampling matrix
\[
\tilde{M}_c = \begin{pmatrix}
  c_2 & 0\\ 0 & c_1
\end{pmatrix}.
\]

The anisotropic sampling of the position parameter is chosen to allow the flexibility to oversample only in the
shearing direction, which will reduce the redundancy of the generated system. This will later become apparent in
the ratio of the associated frame bounds.

Also the range of values for $k$ deserves a comment, since the two values $-2^{j/2}-1$ and $2^{j/2}+1$
which arose from the sampling procedure were replaced by $-\lceil{2^{j/2}}\rceil$ and $\lceil{2^{j/2}}\rceil$.
The reason for this is that $2^{j/2}$ is not an integer for odd scales $j$. Also, for the classical almost-separable
shearlet generator introduced in Example \ref{exa:classicalshearlet}, the
support of the Fourier transform of the associated shearlets does
not anymore fall into $\cC_1 \cup \cC_3$ for $k = -2^{j/2}-1$ and $k = 2^{j/2}+1$ (if considering the horizontal cone).

This now gives the following discrete system.

\begin{definition}
\label{defi:discreteshearlets}
For some sampling vector $c = (c_1,c_2) \in (\RR^+)^2$, the {\em (cone-adapted) regular discrete shearlet system}
$\cSH(c;\phi,\psi,\tilde{\psi})$ generated by a {\em scaling function} $\phi \in L^2(\mathbb{R}^2)$ and
{\em shearlets} $\psi, \tilde{\psi} \in L^2(\mathbb{R}^2)$ is defined by
\[
\cSH(c;\phi,\psi,\tilde{\psi}) = \Phi(c_1,\phi) \cup \Psi(c,\psi) \cup \tilde{\Psi}(c,\tilde{\psi}),
\]
where
\[
\Phi(c_1,\phi) = \{\phi_m  = \phi(\cdot-c_1m) : m \in \bZ^2\},
\]
\[
\Psi(c,\psi) = \{\psi_{j,k,m} =  2^{3j/4} {\psi}({S}_{-k} {A}_{2^j}\cdot-M_c m) :
j \ge 0, |k| \le \lceil 2^{j/2} \rceil, m \in \bZ^2 \},
\]
and
\[
\tilde{\Psi}(c,\tilde{\psi}) = \{\tilde{\psi}_{j,k,m} =  2^{3j/4} \tilde{\psi}(S^T_{-k} \tilde{A}_{2^j}\cdot-\tilde{M}_c m) :
j \ge 0, |k| \le \lceil 2^{j/2} \rceil, m \in \bZ^2 \}.
\]
Setting
\[
\Lambda_{cone}=\{(j,k,m) : j \ge 0, |k| \le \lceil 2^{j/2} \rceil, m \in \bZ^2\},
\]
the associated {\em (Cone-Adapted) Regular Discrete Shearlet Transform} $\mathcal{S}\mathcal{H}_{\phi,\psi,\tilde{\psi}} f : \ZZ^2 \times \Lambda_{cone}^2 \to \mathbb{C}^3$
of some function $f \in L^2(\RR^2)$ is given by
\[
\mathcal{S}\mathcal{H}_{\phi,\psi,\tilde{\psi}} f(m',(j,k,m),(\tilde{j},\tilde{k},\tilde{m}))
 = (\langle f,\phi_{m'}\rangle,\langle f,\psi_{j,k,m}\rangle,\langle f,\tilde{\psi}_{\tilde{j},\tilde{k},\tilde{m}}\rangle).
\]
\end{definition}
The reader should keep in mind that although not indicated by the notation, the functions $\phi_m$, $\psi_{j,k,m}$, and
$\tilde{\psi}_{j,k,m}$ all depend on the sampling constants $c_1, c_2$.
Notice further that the particular sampling set \eqref{eq:regularsampling} forces a change in the ordering of parabolic
scaling and shearing, which will be mimicked later by the class of irregular parameters we will be considering. The
definition itself shows that the shearlets $\psi_{j,k,m}$ (also $\tilde{\psi}_{j,k,m}$) live on anisotropic regions
of size $2^{-j} \times 2^{-\frac j2}$ at various orientations in spatial domain and hence on
anisotropic regions of size $2^{j} \times 2^{\frac j2}$ in frequency domain (see also Figure \ref{fig:tilings}(b)).

\begin{example}
\label{exa:classicalshearlet}
The {\em classical example} of a generating shearlet is a function $\psi \in L^2(\mathbb{R}^2)$
satisfying
\[
\hat{\psi}(\xi) = \hat{\psi}(\xi_1,\xi_2) = \hat{\psi}_1(\xi_1) \, \hat{\psi}_2(\tfrac{\xi_2}{\xi_1}),
\]
where $\psi_1 \in L^2(\bR)$ is a discrete wavelet, i.e., satisfies the discrete Calder\'{o}n condition
given by
$\sum_{j \in \bZ}|\hat\psi_1(2^{-j}\omega)|^2 = 1$ for a.e. $\omega \in \bR$, with
$\hat{\psi}_1 \in C^\infty(\mathbb{R})$ and supp $\hat{\psi}_1 \subseteq [-\frac54,-\frac14] \cup [\frac14,\frac54]$,
and $\psi_2 \in L^2(\mathbb{R})$ is a `bump' function, namely $\sum_{k = -1}^{1}|\hat\psi_2(\omega+k)|^2 = 1$
for a.e. $\omega \in [-1,1]$, satisfying $\hat{\psi}_2 \in C^\infty(\mathbb{R})$
and $\text{\em supp}\,\hat{\psi}_2 \subseteq [-1,1]$. There are several choices of $\psi_{1}$ and $\psi_{2}$
satisfying those conditions, and we refer to \cite{GKL06} for further details.

The tiling of frequency domain given by this band-limited generator and choosing $\tilde{\psi}=\psi(x_2,x_1)$
is illustrated in Figure \ref{fig:tilings}(b). For this particular choice, using an appropriate scaling
function $\phi$ for $\cR$, it was proven in \cite[Thm. 3]{GKL06} that the associated (cone-adapted) discrete
regular shearlet system $\cSH(1;\phi,\psi,\tilde{\psi})$ forms a Parseval frame for $L^2(\RR^2)$.
\end{example}


\subsection{Irregular Sampling of Shearlet Parameters}
\label{subsec:irregularsampling}

Questions arising from sampling theory and the inevitability of perturbations yield the
necessity to study shearlet systems with arbitrary discrete sets of parameters. To derive a better
understanding of this problem, we take the viewpoint of regarding the previously
considered set of parameters
\[
\{(2^{-j},k2^{-j/2},S_{k2^{-j/2}}A_{2^{-j}}M_c m
: j \ge 0,\, k \in \{-\lceil 2^{\frac j2} \rceil,\ldots, \lceil 2^{\frac j2} \rceil \},\, m \in \bZ^2\}
\]
as a discrete subset of $\bS_{cone}$ (for $\cC_1 \cup \cC_3$ and similarly for $\cC_1 \cup \cC_3$ with
$M_c$ substituted by $\tilde{M}_c$) and for the low frequency part $\RR^2$ likewise.
Then, an arbitrary set of parameters is merely a different sampling set, and questions
of density conditions of such a sampling set or its geometric properties are lurking in the
background. Some density-like properties will indeed come into play in this paper mixed with decay
conditions on the shearlet generator. For an extensive analysis of the geometric properties of the
sampling set yielding necessary conditions for shearlet frames directly arising from a group
representation, we refer to \cite{KKL10}.

Next we formally define (cone-adapted) irregular shearlet systems.
\begin{definition}
Let $\Delta$ and $\Lambda$, $\tilde{\Lambda}$ be discrete subsets of $\RR^2$ and $\bS_{cone}$, respectively, and let
$\phi \in L^2(\mathbb{R}^2)$ as well as $\psi, \tilde{\psi} \in L^2(\mathbb{R}^2)$. Then the
{\em (cone-adapted) irregular discrete shearlet system} $\cSH(\Delta,\Lambda,\tilde{\Lambda};\phi,\psi,\tilde{\psi})$ is
defined by
\[
\cSH(\Delta,\Lambda,\tilde{\Lambda};\phi,\psi,\tilde{\psi}) = \Phi(\Delta,\phi) \cup \Psi(\Lambda,\psi) \cup
\tilde{\Psi}(\tilde{\Lambda},\tilde{\psi}),
\]
where
\[
\Phi(\Delta,\phi) = \{\phi_t  = \phi(\cdot-t) : t \in \Delta\},
\]
\[
\Psi(\Lambda,\psi) = \{\psi_{a,s,t}=  a^{-\frac34} \psi(A_a^{-1}S_s^{-1}(\, \cdot \, -t))
: (a,s,t) \in \Lambda\},
\]
and
\[
\tilde{\Psi}(\tilde{\Lambda},\tilde{\psi}) = \{\tilde{\psi}_{a,s,t}=  a^{-\frac34} \tilde{\psi}(\tilde{A}_a^{-1}S_s^{-T}(\, \cdot \, -t))
: (a,s,t) \in \tilde{\Lambda}\}.
\]
Then the associated {\em (Cone-Adapted) Irregular Discrete Shearlet Transform} $\mathcal{S}\mathcal{H}_{\phi,\psi,\tilde{\psi}} f : \Delta \times \Lambda \times \tilde{\Lambda} \to \mathbb{C}^3$
of some function $f \in L^2(\RR^2)$ is given by
\[
\mathcal{S}\mathcal{H}_{\phi,\psi,\tilde{\psi}} f(t',(a,s,t),(\tilde{a},\tilde{s},\tilde{t}))
= (\langle f,\phi_{t'}\rangle,\langle f,\psi_{a,s,t}\rangle,\langle f,\tilde{\psi}_{\tilde{a},\tilde{s},\tilde{t}}\rangle).
\]
\end{definition}


\subsection{Class of Irregular Shearlet Systems}
\label{subsec:class}

Since the low frequency part is already very well studied, and since both cone pairs
$\cC_1$ and $\cC_3$ as well as $\cC_2$ and $\cC_4$ are treated similarly, we from now on
restrict our attention to the horizontal cone
\[
\cC:=\cC_1 \cup \cC_3,
\]
and are interested in frame properties of the system $\Psi(\Lambda,\psi)$, when expanding
functions in
\[
L^2(\cC) = \{f \in L^2(\bR^2) : \supp \wq{\hat{f}} \subseteq \cC\}.
\]
We remark that it is sufficient to consider the function space $L^2(\cC)$, since in
classical shearlet theory, the representing
shearlets are orthogonally projected onto the cones (and, in particular, onto the cone $\cC$
under consideration).
For clarity, we further introduce the notion
\[
\cSH(\Lambda,\psi) := \Psi(\Lambda,\psi).
\]

We will now focus on a class of sets of parameters satisfying some weak conditions, which
in a sense endow them with -- loosely speaking -- similar properties as \eqref{eq:regularsampling}.
We remark that, since one main motivation for considering irregular sets of parameters are
stability questions, such a constraint seems very natural. Also, we
impose weak decay conditions on the shearlet generator, which still
include spatially compactly supported shearlets; our main objective.

\subsubsection{Feasible Sets of Parameters}

We now first discuss our hypotheses on the set of parameters. For this, let
\[
\{(a_j,s_{j,k},t_{j,k,m}^c) : j \ge 0, k \in K_j, m \in \bZ^2\},\qquad K_j \subseteq \bZ, \, c=(c_1,c_2) \in (\RR^+)^2,
\]
be an arbitrary discrete set of parameters in $\bS_{cone}$.

For the decreasing sequence $\{a_j\}_{j\ge 0} \subset
\mathbb{R}^{+}$ to act as a scaling parameter, we require a certain
growth condition, which we phrase as follows: For each $\mu \in (0,
1)$, there shall exist some positive integer $p$ such that
\beq \label{eq:aj}
\frac{a_{j+p}}{a_{j}} < \mu \qquad \textrm{for all } j \ge 0.
\eeq
Obviously, this condition is satisfied for the customarily utilized sequence $a_j = 2^{-j}$ with $p=1$ and
$\mu = 1/2$.

Further, we assume the sequence $\{s_{j,k}\}_{j,k}\subset \mathbb{R}$ of
shear parameters to be of the form
\[
s_{j,k} = s_k \sqrt{a_j}\qquad \mbox{for each }j \ge 0, \, k \in K_j,
\]
which allows a change of the parabolic scaling and shear operators similar to
the discretization of the (cone-adapted) continuous shearlet systems
(see Definition \ref{defi:discreteshearlets}). We also restrict to
\beq\label{eq:sk0}
|s_k| \le a_j^{1/2}+1\qquad \mbox{for all }j \ge 0, \, k \in K_j,
\eeq
and mimic the discretization in the regular case. From now on, we assume that
$0 \in K_j$ for each $j \ge 0$ and $s_0 = 0$ for simplicity, and remark that all of the results in this paper can be
easily extended to a general sequence of shear parameters $s_k$ without this assumption.
Finally, we assume the sequence $\{s_k\}_{k}$ to satisfy a density-like
condition, which will ensure that the shearing operation does not cause a `piling up' of
the essential supports of the shearlets. In particular, we require that
\begin{equation} \label{eq:sk2}
\sup\limits_{x=(x_1,x_2) \in \RR^2} \sum\limits_{k \in \ZZ}  \min\{1,|(S_{s_k}x)_2|\}
\cdot \min{\{1,|(S_{s_k}x)_1|^{-\gamma}\}} \leq C(\gamma) < \infty \qquad\textrm{for any } \gamma >1.
\end{equation}

The sequence $\{t_{j,k,m}^c\}_{j,k,m}\subset \mathbb{R}^2$ of translation parameters shall
have a grid structure, which we impose by assuming that
\beq\label{eq:tjkm}
t_{j,k,m}^c = S_{s_{j,k}}A_{a_j}M_c m,\qquad \mbox{for some } c=(c_1,c_2) \in (\RR^+)^2.
\eeq
Summarizing, we get the following

\begin{definition}
Let the sequences $\{a_j\}_{j \ge 0}\subset \bR^+$, $\{s_{j,k}\}_{j \ge 0, k \in K_j}\subset \mathbb{R}$, and the sequence
$\{t_{j,k,m}^c\}_{j \ge 0, k \in K_j, m \in \bZ^2}\subset \mathbb{R}^2$ satisfy \eqref{eq:aj}, \eqref{eq:sk0}--\eqref{eq:sk2},
and \eqref{eq:tjkm}, respectively. Then we call the sequences
\[
\{(a_j,s_{j,k}= s_k \sqrt{a_j},t_{j,k,m}^c= S_{s_{j,k}}A_{a_j}M_c m) : j \ge 0, k \in K_j, m \in \bZ^2\}
\]
a {\em feasible} set of parameters.
\end{definition}


\subsubsection{Feasibility Conditions for Generating Shearlet}

We next impose the following decay condition on the generating shearlets $\psi \in L^2(\bR^2)$:

\begin{definition}
A function $\psi \in L^2(\RR^2)$ is called a {\em feasible shearlet}, if there exist
$\alpha > \gamma > 3$, and $q > q' >0$,  $q > r >0$ such that
\begin{equation}\label{eq:psi}
|\hat{\psi}(\xi_{1}, \xi_{2})| \le \min\{1,|q\xi_1|^{\alpha}\} \cdot \min{\{1,|q'\xi_1|^{-\gamma}\}} \cdot \min{\{1,|r\xi_2|^{-\gamma}\}}.
\end{equation}
\end{definition}
At first sight it seems picky to choose three different constants $q, q'$, and $r$ for the different
`decay areas' in frequency domain. However, this flexibility will become crucial for deriving
estimates for the lower
frame bound, and hence for analyzing the class of compactly supported shearlet frames we aim to
introduce.

\begin{figure}[ht]
\begin{center}
\hspace*{10pt}
\includegraphics[width=9cm]{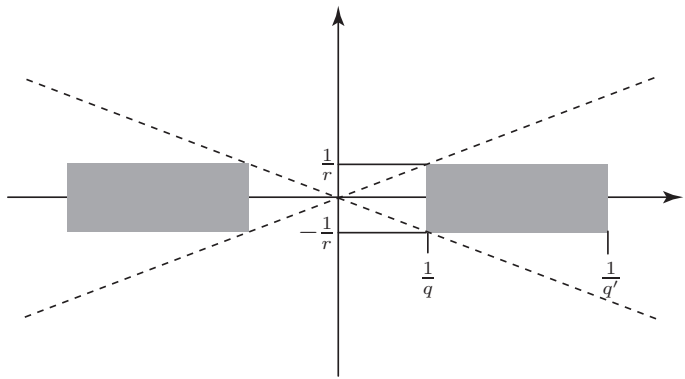}
\put(-100,35){\footnotesize{$\frac{1}{q}$}}
\put(-32,35){\footnotesize{$\frac{1}{q'}$}}
\put(-146,52){\footnotesize{$-\frac{1}{r}$}}
\put(-138,78){\footnotesize{$\frac{1}{r}$}}
\end{center}
\caption{Essential support of $|\hat \psi|$ in the frequency domain for a shearlet $\psi$ satisfying the
decay condition \eqref{eq:psi}.}
\label{fig:decay}
\end{figure}

The class of feasible shearlets certainly includes the classical example of generating shearlets
(Example \ref{exa:classicalshearlet}) as well as all $L^2$-functions whose Fourier transform
is of polynomial decay of a certain degree. In particular, it includes functions
which are spatially compactly supported.
Moreover, condition \eqref{eq:psi} is chosen in such a way that $\psi$ exhibits an essential support
in frequency domain somehow similar to the support of the classical shearlets defined in Example
\ref{exa:classicalshearlet}, but is adapted to separable shearlet frames which later on will
serve as the structure for our explicitly constructed compactly supported shearlet frames.
For an illustration we refer to Figure \ref{fig:decay}.


\subsubsection{Class of Shearlet Systems}

Concluding, in this paper we are concerned with a particular class of (cone-adapted) irregular shearlet
systems defined as follows.

\begin{definition}
\label{defi:feasibility}
Let $\Lambda = \{(a_j,s_{j,k},t_{j,k,m}^c) : j \ge 0, k \in K_j,  m \in \bZ^2\}$ be a feasible set of
parameters and let $\psi \in L^2(\bR^2)$ be a feasible shearlet. Then we call
\begin{eqnarray*}
\mathcal{S}\mathcal{H}(\Lambda,\psi) = \{\psi_{j,k,m} =
a_j^{-\frac34} \psi(S_{-s_k} A_{a_j^{-1}}  \, \cdot \, -M_c m) :  j \ge 0, k \in K_j, m \in \bZ^2\}
\end{eqnarray*}
a {\em feasible shearlet system}.
\end{definition}
Notice that such a system depends in particular on the parameters $\mu, \alpha, \gamma,q , q',$ and $r$
as well as on the function $C(\cdot)$.
We further remark the abuse of notation for $\psi_{j,k,m}$ as compared to Definition \ref{defi:discreteshearlets} to
keep the notation simple. However, it will always be made clear which interpretation is meant.


\section{A General Sufficient Condition for Shearlet Frames}
\label{sec:generalsufficient}

In this section, we will state a very general sufficient condition for the existence of cone-adapted irregular
shearlet frames. A similar condition can be found in \cite{KL07}, however for shearlet frames
arising from a group representation. This result can not be directly carried over, since the
relation between the cone, the bound $a_j^{-1/2}$ for the parameters $s_k$, and the most likely
not band-limited shearlet generator $\psi$ needs careful handling. Also, the result in \cite{KL07}
does not include compactly supported shearlet generators.

\subsection{Covering Properties}

For stating and proving the result, we require a particular function, which is an adapted
form already appeared in \cite{KL07} (see also \cite{Dau92} for the wavelet situation). It
can be regarded as a variant of the main term in the $t_q$ condition for tight affine frames
\cite{Dau92}.

Let $\mathcal{S}\mathcal{H}(\Lambda,\psi)$ be a feasible shearlet system as defined in Definition
\ref{defi:feasibility}.
The main objective will now be the function $\Phi : \cC \times \bR^2 \to \bR$ defined by
\begin{equation}\label{eq:Phi}
\Phi(\xi,\omega) = \sum\limits_{j \ge 0} \sum\limits_{k \in K_j}
|\hat{\psi}(S_{s_{k}}^{T} A_{a_{j}}\xi)| |\hat{\psi}(S_{s_{k}}^{T}A_{a_{j}}\xi + \omega)|.
\end{equation}
Notice that this function measures the extent to which the essential supports of the
scaled and sheared versions of the shearlet generator overlaps. For later use, we
also introduce the function $\Gamma : \bR^2 \to \bR$ defined by
\[
\Gamma(\omega) = \esssup\limits_{\xi \in \cC} \Phi(\xi,\omega),
\]
measuring the maximal extent to which these versions overlap for a fixed distance,
as well as the values
\beq \label{eq:defiLinLsup}
L_{inf} = \essinf\limits_{\xi \in \cC} \Phi(\xi,0)
\qquad \mbox{and} \qquad
L_{sup} = \esssup\limits_{\xi \in \cC} \Phi(\xi,0),
\eeq
which relate to the classical discrete Calder\'{o}n condition. Finally, we also require
the values
\begin{equation} \label{eq:defR}
R(c) = \sum\limits_{m \in \Z^{2}\setminus\{0\}} \left[ \Gamma\left({M_c^{-1}{m}}\right)
\Gamma\left({-M_c^{-1}{m}}\right)\right]^{1/2}, \quad \text{where} \quad c = (c_1,c_2) \in (\RR^+)^2.
\end{equation}
The function $R(c)$ measures the average of the symmetrized function values $\Gamma(M_c^{-1}m)$.

We now first turn our attention to the terms $L_{sup}$ and $R(c)$ and provide explicit estimates for
those, which will later be used for estimates for frame bounds associated to a shearlet system.
We start by estimating $L_{sup}$ in

\begin{proposition} \label{prop:Lsupfinite}
Let $\mathcal{S}\mathcal{H}(\Lambda,\psi)$  be a feasible shearlet system, and let $L_{sup}$ be defined
as in \eqref{eq:defiLinLsup}. Then
\begin{equation}\label{eq:Upp}
L_{sup} \le p \cdot \Bigl( \frac qr \cdot C(2\gamma) \Bigr)
\Bigl(\Big\lceil \log_{1/\mu}\Big(\frac{q}{q'}\Big) \Big\rceil +\frac{1}{1-\mu^{2\alpha-1}}+\frac{1}{1-\mu^{2\gamma}}\Bigr) < \infty.
\end{equation}
In the special case $s_k = k$, $k \in \ZZ$, the following stronger estimate holds:
\begin{equation}\label{eq:Upp2}
L_{sup} \le p \cdot \Bigl(  \frac{q}{r}\left( 2+\frac{2}{2\gamma-1}\right)+1 \Bigr)
\Bigl(\Big\lceil \log_{1/\mu}\Big(\frac{q}{q'}\Big) \Big\rceil +\frac{1}{1-\mu^{2\alpha-1}}+\frac{1}{1-\mu^{2\gamma}}\Bigr) < \infty.
\end{equation}
\end{proposition}

One essential ingredient in the proof of this lemma is a result from \cite{YZ05}, which
we state here for the convenience of the reader.

\begin{lemma} \label{lemma:estimate_aj} \cite[Lem. 2.1]{YZ05}
Let $\mu \in (0, 1), p \in \mathbb{Z}^{+}$, and let $\{a_{j}\}_{j \in \mathbb{Z}} \subset \R^{+}$ be a
decreasing sequence which satisfies (\ref{eq:aj}). Then, for all $\alpha > 0$ and $t > 0$,
\[
\sum\limits_{a_{j} \ge t} a_{j}^{-\alpha} \le p\, t^{-\alpha} \frac{1}{1 - \mu^{\alpha}},
\qquad \sum\limits_{a_{j} \le t} a_{j}^{\alpha} \le p\, t^{\alpha}
\frac{1}{1 - \mu^{\alpha}}.
\]
\end{lemma}

This now allows us to prove Proposition \ref{prop:Lsupfinite}.

\begin{proof}[Proof of Proposition \ref{prop:Lsupfinite}]
Let $\xi = (\xi_{1}, \xi_{2})^{T} \in {\mathbb{R}^{2}}$.
By \eqref{eq:psi}, we obtain
\begin{eqnarray}\nonumber
\lefteqn{\Phi(\xi,0)}\\ \label{eq:estimate1}
& \leq & \hspace*{-0.25cm} \sup_{\xi \in \RR^2} \sum_{j \ge 0}  \min{\{1,|qa_j\xi_1|^{2\alpha}\}}\min\{1,|q'a_j\xi_1|^{-2\gamma}\}
\sum_{k}\min\{1,|r(a_j^{1/2}\xi_2+s_ka_j\xi_1)|^{-2\gamma}\}.
\end{eqnarray}
Letting $\eta_1 = q\xi_1$,
\begin{eqnarray}\nonumber
{\Phi(\xi,0)} &\le& \sup_{(\eta_1,\xi_2) \in \RR^2} \sum_{j \ge 0} \min{\{1,|a_j\eta_1|^{2\alpha-1}\}}
\min\{1,\left|q'q^{-1}a_j\eta_1 \right|^{-2\gamma}\} \sum_{k \in \ZZ} \frac qr \min\left \{1,\left|rq^{-1}a_j\eta_1\right| \right\}\\  \label{eq:Phi0}
&& \;\cdot \min\{1,|ra_j^{1/2}\xi_2+rq^{-1}a_j\eta_1s_k|^{-2\gamma}\}.
\end{eqnarray}
By \eqref{eq:sk2},
\[
\sum_{k \in \ZZ} \frac qr \min  \{ 1,|rq^{-1}a_j\eta_1|\}
\min\{1,|ra_j^{1/2}\xi_2+rq^{-1}a_j\eta_1s_k|^{-2\gamma}\}  \le \frac{q}{r}C(2\gamma),
\]
Hence we can continue \eqref{eq:Phi0} by
\begin{eqnarray*}
\Phi(\xi,0)
&\leq& \Bigl(\frac qr C(2\gamma) \Bigr)\sup_{\eta_1 \in \RR} \sum_{j \ge 0} \min{\{1,|a_j\eta_1|^{2\alpha-1}\}}
\min\{1,|q'q^{-1}a_j\eta_1|^{-2\gamma}\}  \\
&=&  \Bigl(\frac qr C(2\gamma) \Bigr) \sup_{\eta_1 \in \RR} \Bigl(\sum_{j \ge 0}
\chi_{[0,1)}(|a_j\eta_1|)|a_j\eta_1|^{2\alpha-1} + \chi_{[1,\frac{q}{q'})}(|a_j\eta_1|)\\
& & \hspace*{3cm}+
\chi_{[\frac{q}{q'},\infty)}(|a_j\eta_1|)|q'q^{-1}a_j\eta_1|^{-2\gamma}\Bigl)\\
&\leq&   \Bigl(\frac qr C(2\gamma) \Bigr) \cdot \sup_{\eta_1 \in \RR} \Bigl(
\sum_{|a_j\eta_1| \leq 1}|a_j\eta_1|^{2\alpha-1}+ \sum_{j\ge 0} \chi_{[1,\frac{q}{q'}]}(|a_j\eta_1|)+\hspace*{-0.5cm}
\sum_{|q'q^{-1}a_j\eta_1| \ge 1} |q'q^{-1}a_j\eta_1|^{-2\gamma}\Bigr).
\end{eqnarray*}
The first claim \eqref{eq:Upp} now follows from Lemma \ref{lemma:estimate_aj}.

Next assume that $s_k = k$ for all $k \in \ZZ$. To derive an estimate for $L_{sup}$ in this case, we first
split \eqref{eq:estimate1} into the cases $k=0$ and $k \neq 0$ with $\eta_1 = q\xi_1$ and $\eta_2 = ra_j\xi_2$,
which yields
\begin{eqnarray}\nonumber
\lefteqn{\Phi(\xi,0)} \\ \nonumber
&\le& \sup_{\eta_1 \in \RR} \sum_{j \ge 0} \min{\{1,|a_j\eta_1|^{2\alpha-1}\}}\min\{1,|q'q^{-1}a_j\eta_1|^{-2\gamma}\}
\sup_{\eta_2 \in [0,|rq^{-1}a_j\eta_1|]}\sum_{k \neq 0} \frac qr \min\{1,|rq^{-1}a_j\eta_1|\} \\  \label{eq:estimate2}
&& \;\cdot \min\{1,|\eta_2+rq^{-1}a_j\eta_1 k|^{-2\gamma}\} + \sup_{\eta_1 \in \RR} \sum_{j \ge 0} \min{\{1,|a_j\eta_1|^{2\alpha}\}}
\min\{1,|q'q^{-1}a_j\eta_1|^{-2\gamma}\}.
\end{eqnarray}
Notice that
\[
\sup_{\eta_2 \in [0,|rq^{-1}a_j\eta_1|]}\sum_{k \neq 0} \frac qr \min\{1,|rq^{-1}a_j\eta_1|\}
\min\{1,|\eta_2+rq^{-1}a_j\eta_1 k|^{-2\gamma}\} \le \frac{q}{r}\left( 2+\frac{2}{2\gamma-1}\right)
\]
and
\[
\min{\{1,|a_j\eta_1|^{2\alpha}\}} \min\{1,|q'q^{-1}a_j\eta_1|^{-2\gamma}\}
\le \min{\{1,|a_j\eta_1|^{2\alpha-1}\}}\min\{1,|q'q^{-1} a_j\eta_1|^{-2\gamma}\}.
\]
Therefore, we can continue \eqref{eq:estimate2} by
\[
\Phi(\xi,0)
\leq \Bigl(  \frac{q}{r}\left( 2+\frac{2}{2\gamma-1}\right)+1 \Bigr)\sup_{\eta_1 \in \RR} \sum_{j \ge 0}
\min{\{1,|a_j\xi_1|^{2\alpha-1}\}}\min\{1,|q'q^{-1}a_j\eta_1|^{-2\gamma}\}.
\]
Now \eqref{eq:Upp2} follows from Lemma \ref{lemma:estimate_aj}; the proposition is proved.
\end{proof}

The strengthened estimate in the special case of a shearing parameter sequence $\{s_k\}_k$ with $s_k = k$ for $k \in \ZZ$
will become important when estimating the frame bounds of the concrete class of compactly supported shearlet frames we will
introduce in Section \ref{sec:compactsupport}.

The next result reveals the dependence of $R(c)$ on the sampling matrix $M_c$ and provides a useful
estimate for such. For its very technical proof, we refer to Subsection \ref{subsubsec:estimateR}.

\begin{proposition}
\label{prop:estimateR}
Let $\mathcal{S}\mathcal{H}(\Lambda,\psi)$ be a feasible shearlet system with associated sampling vector $c = (c_1,c_2)
\in (\RR^+)^2$ satisfying $c_1 \ge c_2$, and let $R(c)$ be defined as in \eqref{eq:defR}.
Then, for any $\gamma'$, which satisfies $1<\gamma'<\gamma-2$, we have
\begin{equation}\label{eq:upper_R}
R(c) \leq T_1 D_1(\gamma) + \min \left \{\left\lceil \frac{c_1}{c_2} \right\rceil,2 \right \} T_2 D_2(\gamma - \gamma') + T_3 (D_1(\gamma)+D_2(\gamma)),
\end{equation}
where
\begin{eqnarray*}
T_1 \hspace*{-0.1cm}&\hspace*{-0.1cm}=\hspace*{-0.1cm}&\hspace*{-0.1cm} p\Bigl( \frac{q}{r}C(\gamma)\Bigr)
\Bigl(\Bigl\lceil \log_{\frac{1}{\mu}}\Bigl(\frac{q}{q'}\Bigr)\Bigl\rceil+\frac{1}{1-\mu^{\gamma}}+\frac{1}{1-\mu^{\alpha-\gamma}}\Bigr)
\Bigl( \frac{2c_1}{q'}\Bigr)^{\gamma}, \\
T_2 \hspace*{-0.1cm}&\hspace*{-0.1cm}=\hspace*{-0.1cm}&\hspace*{-0.1cm} p \Bigl( \frac{q}{r}C(\gamma')\Bigr)
\Bigl(2\Bigl\lceil \log_{\frac{1}{\mu}}\Bigl(\frac{q}{q'}\Bigr)\Bigl\rceil+\frac{1}{1-\mu^{\gamma'}}+\frac{1}{1-\mu^{\alpha-\gamma'}}
+
\frac{1}{1-\mu^{\gamma}}+\frac{1}{1-\mu^{\alpha-\gamma}}\Bigr)\Bigl( \frac{2qc_2}{q'r}\Bigr)^{\gamma-\gamma'}, \\
T_3 \hspace*{-0.1cm}&\hspace*{-0.1cm}=\hspace*{-0.1cm}&\hspace*{-0.1cm} p\Bigl( \frac{q}{r}C(\gamma)\Bigr)
\Bigl( \frac{1}{1-\mu^{\gamma}}\Bigr)\Bigl(\frac{2c_1}{q'} \Bigr)^{\gamma},
\end{eqnarray*}
and, for any $h > 0$,
\[
D_1(h) = 2\Bigl(1+\frac{1}{h-1}\Bigl)+\frac{4}{h-1}\Bigl(1+\frac{1}{h-2}\Bigr), \quad
D_2(h) = 6\Bigl(1+\frac{1}{h-1}\Bigl)+\frac{4}{h-1}\Bigl(1+\frac{1}{h-2}\Bigr).
\]
In particular, for any $\gamma'$, which satisfies $1<\gamma'<\gamma-2$, there exist positive constants
$\kappa_1$ and $\kappa_2$ (in particular, independent on $c$) such that
\begin{equation}\label{eq:decayR}
R(c) \le \kappa_1\left( \frac{2c_1}{q'}\right)^{\gamma}+\kappa_2\left(\frac{2qc_2}{q'r}\right)^{\gamma-\gamma'} .
\end{equation}
\end{proposition}


\subsection{General Sufficient Condition}

Retaining the notions introduced in the previous subsection, we can now formulate
a general sufficient condition for the existence of cone-adapted irregular
shearlet frames. We defer the lengthy proof to Subsection \ref{subsubsec:general}.

\begin{theorem}
\label{theo:general}
Let $\mathcal{S}\mathcal{H}(\Lambda,\psi)$ be a feasible shearlet system, let $L_{inf}$, $L_{sup}$ be defined
as in \eqref{eq:defiLinLsup}, and let $R(c)$ be defined as in \eqref{eq:defR}. If $\tilde{L}_{sup}$ and $\tilde{L}_{inf}$
are chosen such that
\[
R(c) < \tilde{L}_{inf} \le L_{inf} \,\, \mbox{and} \,\, L_{sup} \le \tilde{L}_{sup},
\]
then $\mathcal{SH}(\Lambda, \psi)$ is a frame for $L^2(\cC)$ with frame bounds $A$ and $B$ satisfying
\[
\frac{1}{|\det M_c|} [\tilde{L}_{inf} - R(c)] \le A \le B \le \frac{1}{|\det M_c|} [\tilde{L}_{sup} + R(c)].
\]
\end{theorem}

The role of the determinant of the sampling matrix $M_c$, which can also be regarded as the inverse of the Beurling
density of the translation sampling grid in spatial domain, requires a careful examination.
The `quality' of a frame can be measured by the magnitude of the quotient of the frame bounds
$B/A$, since the closer it is to $1$, the better frame reconstruction algorithms perform.
In this case, this quotient is bounded by
\beq \label{eq:estimateAB}
\frac{B}{A} \le \frac{\tilde{L}_{sup} + R(c)}{\tilde{L}_{inf} - R(c)}.
\eeq
We remark that in order to estimate the quotient $B/A$, it is sufficient to estimate
$R(c)$ and $\tilde{L}_{inf}$ in \eqref{eq:estimateAB}, since we can use Proposition \ref{prop:Lsupfinite} for $L_{sup}$.

We make the following observations:
\bitem
\item[(a)] Notice that as $c_1$ approaches $0$ -- allowed by \eqref{eq:decayR} -- the
quotient $B/A$ approaches $\tilde{L}_{sup}/\tilde{L}_{inf}$, which is optimized by choosing $\tilde{L}_{inf}=L_{inf}$ and
$\tilde{L}_{sup} = L_{sup}$. This is
very intuitive, since as $c_1 \to 0$ the translation grid becomes denser, the frame property
more likely, and it becomes merely a question depending on the values $L_{inf}$ and $L_{sup}$.
\item[(b)] \wq{As we indicated in the previous remark, one can achieve better frame bounds by
simply choosing a sufficiently small sampling constant $c_1$ ($= \max\{c_1,c_2\}$). However, this leads to a highly
redundant system, which causes a significant computational complexity for practical applications.
Therefore, we might not want to choose $c_1$ too small to avoid high redundancy,} but instead seek a
threshold, which gives the largest possible value for $c_1$ and also $c_2$.
\item[(c)] Reasoning about the previous two observations it becomes apparent that
the balance between the possible values of $c_1$ and $c_2$ is an optimization problem, which presumably
is different depending on the application at hand.
\eitem

\section{Compactly Supported Shearlet Frames}
\label{sec:compactsupport}

We must now brace ourselves for the considerably harder challenge to construct compactly supported
(separable) shearlet frames with provably `reasonable' frame bounds. However, this task seems, at least to us,
much more rewarding as it shows the delicacies of exploiting the previously derived general
sufficient frame conditions.

The main idea will be to start with a maximally flat low pass filter and use this to build a
separable shearlet generator, which is then shown to generate a shearlet frame with `good' frame
bounds -- in the sense of a small ratio. The modeling of the shearlet mimics the classical Example
\ref{exa:classicalshearlet} with a wavelet-like function in horizontal direction and a bump-like function
in vertical direction (on $\CC$), while ensuring spatial compact support.
This is first performed for $L^2(\cC)$, followed by a discussion of $L^2(\RR^2)$.
The new shearlet frame is finally shown to provide (almost) optimally sparse approximations for
cartoon-like images.

\subsection{Shearlet Frame for $L^2({\cC})$}
\label{subsec:compact}

We start by constructing a compactly supported (separable) shearlet frame $\mathcal{S}\mathcal{H}(\Lambda,\psi)$
for $L^2(\cC)$. Our construction will ensure that $\mathcal{S}\mathcal{H}(\Lambda,\psi)$ is a feasible shearlet
system, which enables us to exploit our previous results, in particular, Proposition \ref{prop:Lsupfinite}, Theorem \ref{theo:general},
and Proposition \ref{prop:estimateR} to derive estimates of its frame bounds. For our construction, we
focus on regular discrete shearlet systems, i.e., we choose
\[
\Lambda = \{(2^{-j},k2^{j/2},S_{k2^{j/2}}A_{2^{j}}M_c m) : j \ge 0,|k| \le \lceil 2^{j/2} \rceil, m \in \bZ^2\}
\qquad \mbox{with } c=(c_1,c_2) \in (\RR^+)^2,
\]
hence the main objective will be to construct a compactly supported shearlet $\psi$, which satisfies the
decay condition \eqref{eq:psi}.

\begin{figure}[ht]
\begin{center}
\includegraphics[height=2in, width=3.5in]{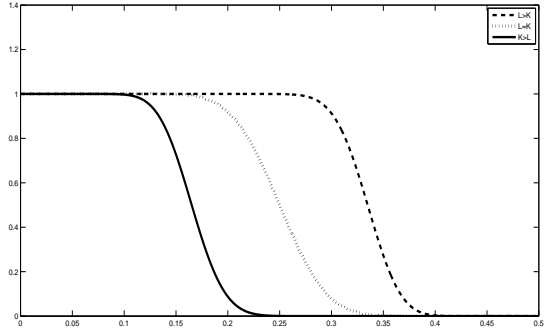}
\end{center}
\caption{Graph of $|m_0|^2$ with $K>L$ (solid), $K=L$ (dot), and $K<L$ (dash).}
\label{fig:m0}
\end{figure}

Our first goal towards this aim is to define a particular univariate scaling function, which will become a main building block
for the compactly supported shearlet $\psi$ we aim to construct, and prove controllable estimates for its function
values. For this, given $K, L \in \bZ^+$, we first let $m_0$ be the trigonometric polynomial with real
coefficients which satisfies $m_0(0) = 1$ and
\begin{equation}\label{eq:m0}
|m_0(\xi_1)|^2 = (\cos(\pi\xi_1))^{2K} \sum_{n=0}^{L-1}{K-1+n\choose n}(\sin(\pi\xi_1))^{2n}.
\end{equation}
Notice that $m_0$ can be obtained from \eqref{eq:m0} using spectral factorization (see \cite{Dau92}).
The two parameters $K$ and $L$ determine the level
of flatness near $\frac12$ and $0$, respectively. In particular, if $K=L$, then $|m_0(\xi_1)|^2$ is symmetric
with respect to the point $\xi_1 = \frac14$, and in this case $m_0$ generates the Daubechies scaling function \cite{Dau92}.
The more general trigonometric polynomial $m_0$ we are considering here is sometimes also referred to
as the {\em maximally flat low pass filter}. For an illustration of $|m_0|^2$ we refer to
Figure \ref{fig:m0}.

The following result states some useful properties, which will become handy in the sequel. For the lengthy proof, we
refer to Subsection \ref{subsec:proofm0}.

\begin{lemma}\label{lemma:m0}
Retaining the previously introduced notations, we have the following conditions.
\bitem
\item[(i)] $|m_0|^2$ is even.
\item[(ii)] $|m_0|^2$ is decreasing on $(0,\frac12).$
\item[(iii)] If $\frac{L-1}{K+L-2} \ge \frac14$, then $|m_0|^2$ is concave on $(0,\frac16)$.
\eitem
\end{lemma}

The next step is to define the sought scaling function $\phi \in L^2(\RR)$ by
\begin{equation}\label{eq:phi}
\hat{\phi}(\xi_1) = \prod_{j=0}^{\infty} m_0(2^{-j}\xi_1).
\end{equation}
Our first result investigating this function concerns a lower bound for the function
$|\hat{\phi}|^2$ on $[-\frac16,\frac16]$.

\begin{proposition}\label{proposition:piLower}
Retaining the previously introduced notations, if $\frac{L-1}{K+L-2} \ge \frac14$, then
\beq \label{eq:lowerphi}
|\hat{\phi}(\xi_1)|^2 \ge \prod_{j=0}^{J_0-1}\left|m_0\left(\tfrac{2^{-j}}{6}\right)\right|^2 \cdot e^{-2^{-J_0+2}(1-|m_0(\frac16)|^2)}
\cdot \chi_{[-1/6,1/6]}(\xi_1), \qquad \xi_1 \in \RR,
\eeq
where $J_0$ is chosen sufficiently large. Moreover, the function $\phi$ is compactly supported.
\end{proposition}

\begin{proof}
First note, that it suffices to consider the interval $[0,\frac16]$, since $|\hat{\phi}(\xi_1)|^2$ is even
(see Lemma \ref{lemma:m0}(i)). Now assume that $\xi_1 \in [0,\frac16]$.
Since, by Lemma \ref{lemma:m0}(ii)+(iii), $|m_0|^2$ is monotone and concave on $(0,\frac16)$, we have
\[
|m_0(\xi_1)|^2 \ge |m_0(0)|^2-\left[6(|m_0(0)|^2-|m_0(\tfrac16)|^2)\right]|\xi_1|.
\]
Since by construction $|m_0(0)|=1$, this implies
\[
|m_0(\xi_1)|^2 \ge 1-6C_1|\xi_1|,
\]
where
\[
C_1 = 1-|m_0(\tfrac16)|^2>0.
\]
Let now $J_0>0$ be chosen sufficiently large such that $2^{-J_0}C_1 \le \frac 12$. Then, since
$1-x \ge e^{-2x}$ for $x \in [0,\frac12]$, we can continue our estimation by
\[
|m_0(2^{-J_0}\xi_1)|^2 \ge e^{-12(2^{-J_0}C_1)|\xi_1|}.
\]
By again using Lemma \ref{lemma:m0}(ii), this finally implies
{\allowdisplaybreaks
\begin{eqnarray*}
|\hat{\phi}(\xi_1)|^2
&=& \prod_{j=0}^{J_0-1}|m_0(2^{-j}\xi_1)|^2\prod_{j'=0}^{\infty}|m_0(2^{-J_0-j'}\xi_1)|^2\\
&\ge& \prod_{j=0}^{J_0-1}\left|m_0\left(\tfrac{2^{-j}}{6}\right)\right|^2\prod_{j'=0}^{\infty}e^{(-12C_1)(2^{-j'-J_0}\xi_1)} \\
&=& \prod_{j=0}^{J_0-1}\left|m_0\left(\tfrac{2^{-j}}{6}\right)\right|^2e^{-2^{-J_0+2}(6C_1)\xi_1} \\
&\ge& \prod_{j=0}^{J_0-1}\left|m_0\left(\tfrac{2^{-j}}{6}\right)\right|^2e^{-2^{-J_0+2}C_1}.
\end{eqnarray*}
}

Finally, by standard arguments (cf. \cite{Dau92}), the function $\phi$ is compactly supported.
The proposition is proved.
\end{proof}

In addition, an upper, `mathematically controllable' bound for the values of $|\hat{\phi}|^2$ will be
required for the analysis of the sought compactly supported shearlet $\psi$. We start with some preparation for being able
to state the precise estimate. For this, we need to introduce a trigonometric polynomial $\tilde{m}_0$
slightly differing from the previously discussed $m_0$. Letting $K' \in \bZ$ be such that $0<K'<K$,
we set
\[
|\tilde{m}_0(\xi_1)|^2 = \left(\cos(\pi\xi_1)\right)^{2K'}\sum_{n=0}^{L-1}{K-1+n\choose n}(\sin(\pi \xi_1))^{2n}
\]
and again can obtain $\tilde{m}_0$ from this using spectral factorization.
The next estimate for the values of $|\tilde{m}_0|^2$ will be crucial for our upper estimate for
$|\hat{\phi}|^2$. The very technical proof is deferred to Subsection \ref{subsec:prooftildem}.

\begin{lemma}\label{lemma:tildem}
Assume that $L \ge 6$ and $L+1 \le K \le 3L-2$ in the definition of $m_0$, i.e., in \eqref{eq:m0}.
Let $K'$ be a non-negative integer which satisfies $\frac{K-K'}{K'+L-1} \ge \frac 14$.
Retaining the previously introduced notations, we have the following conditions.
\bitem
\item[(i)] $|\tilde{m}_0|^2$ is increasing on $(0,\frac16)$.
\item[(ii)] We have
\beq \label{eq:c2}
\max_{\xi_1 \in [0,1]}|\tilde{m}_0(\xi_1)|^2 \le
\sum_{n=0}^{L-1}{K-1+n\choose n}\Bigl( \frac{K'}{K'+n}\Bigr)^{K'}\Bigl(\frac{n}{K'+n}\Bigr)^{n} =: C_2.
\eeq
\eitem
\end{lemma}

This now enables us to state the upper estimate for $|\tilde{m}_0|^2$  we aimed for.

\begin{proposition}\label{proposition:piUpper}
Retaining the previously introduced notations, and letting the non-negative integers $K, K',$ and $L$ be chosen such that $L \ge 6$,
$L+1 \le K \le 3L-2$, and $\frac{K-K'}{K'+L-1} \ge \frac 14$, then, for $J_1$ sufficiently large,
\[
|\hat{\phi}(\xi_1)|^2 \le \min\left\{1,C_2 \cdot |2\pi\xi_1|^{-2\gamma}
\cdot\prod_{j=0}^{J_1-1}\left|\tilde{m}_0\left(\frac{2^{-j}}{2\pi}\right)\right|^2
\cdot e^{2^{-J_1+1}\sum_{n} |h(n)||n|}\right\}, \qquad \xi_1 \in \RR,
\]
where $\gamma = (K-K')-\frac{1}{2}\log_2(C_2)$ and $h(n)$ being the Fourier coefficients of $|\tilde{m}_0|^2$.
\end{proposition}

\begin{proof}
First, observe that $|\hat{\phi}(\xi_1)|^2 \leq 1$ for all $\xi_1 \in \bR$, since $|m_0(\xi_1)|^2 \le 1$.
Now fix $\xi_1 \in \RR$.

If $|2\pi\xi_1| \le 1$, we use the easy conclusion from the proof of Lemma \ref{lemma:tildem}(i),
\[
C_2 \cdot \prod_{j=0}^{J_1-1}\left|\tilde{m}_0\left(\frac{2^{-j}}{2\pi}\right)\right|^2
e^{2^{-J_1+1}\sum_{n} |h(n)||n|} \ge 1,
\]
to show that
\[
\min\left\{1,C_2 \cdot |2\pi\xi_1|^{-2\gamma}
\cdot\prod_{j=0}^{J_1-1}\left|\tilde{m}_0\left(\frac{2^{-j}}{2\pi}\right)\right|^2
\cdot e^{2^{-J_1+1}\sum_{n} |h(n)||n|}\right\} = 1 \ge |\hat{\phi}(\xi_1)|^2,
\]
which was claimed.

Now assume that $|2\pi\xi_1| \ge 1$. First, we observe that
\[
|m_0(\xi_1)|^2 = {\cos(\pi\xi_1)}^{2(K-K')}|\tilde{m}_0(\xi_1)|^2
\]
and
\[
|\hat{\phi}(\xi_1)|^2 = \prod_{j=0}^{\infty}\left( \cos(\pi2^{-j}\xi_1)\right)^{2(K-K')}|\tilde{m}_0(2^{-j}\xi_1)|^2.
\]
Using the classical formula $\prod_{j=1}^{\infty}\cos(2^{-j}x) = \frac{\sin(x)}{x}$ (see \cite{Dau92}), it
follows that
\begin{equation}\label{eq:phi_sin}
|\hat{\phi}(\xi_1)|^2 = \Bigl| \frac{\sin(2\pi\xi_1)}{2\pi\xi_1}\Bigr|^{2(K-K')}\prod_{j=0}^{\infty}|\tilde{m}_0(2^{-j}\xi_1)|^2.
\end{equation}
Since $|2\pi\xi_1| \ge 1$, there exists a positive integer $J$ such that $2^{J-1} \le 2\pi|\xi_1| \le 2^J$.
Fixing this $J$, by \eqref{eq:phi_sin} and the definition of $C_2$, we have
\begin{eqnarray*}
|\hat{\phi}(\xi_1)|^2
& = & |2\pi\xi_1|^{-2(K-K')}\prod_{j'=0}^{J-2}|\tilde{m}_0(2^{-j'}\xi_1)|^2\prod_{j=J-1}^{\infty}
|\tilde{m}_0(2^{-j}\xi_1)|^2\\
& \le & |2\pi\xi_1|^{-2(K-K')}\prod_{j'=0}^{J-2}2^{\log_2(C_2)}\prod_{j=0}^{\infty}|\tilde{m}_0(2^{-j}(2^{-J+1}\xi_1))|^2\\
&\le& |2\pi\xi_1|^{-2(K-K')}\left(2^{(J-1)\log_2(C_2)}\right)\prod_{j=0}^{\infty}|\tilde{m}_0(2^{-j}(2^{-J+1}\xi_1))|^2.
\end{eqnarray*}
Since $2^{J-1} \le 2\pi|\xi_1| \le 2^J$ implies $2^{(J-1)\log_2(C_2)} \le |2\pi\xi_1|^{\log2(C_2)}$,
letting $\eta_1 = 2^{-J+1}\xi_1$, hence $\frac{1}{2\pi} \le |\eta_1| \le \frac{1}{\pi}$, yields
\begin{eqnarray*}
|\hat{\phi}(\xi_1)|^2
&\le& |2\pi\xi_1|^{-2(K-K')+\log_2(C_2)}\sup_{\eta_1 \in [0,1/\pi]}\prod_{j=0}^{\infty}|\tilde{m}_0(2^{-j}\eta_1)|^2\\
&\leq& C_2 \cdot |2\pi\xi_1|^{-2(K-K')+\log_2(C_2)}\sup_{\eta_1 \in [0,1/2\pi]}\prod_{j=0}^{\infty}|\tilde{m}_0(2^{-j}\eta_1)|^2.
\end{eqnarray*}
By Lemma \ref{lemma:tildem}(i), this implies that, for all $J_1 > 0$,
\[
|\hat{\phi}(\xi_1)|^2 \le C_2\cdot |2\pi\xi_1|^{-2\gamma}\prod_{j=0}^{J_1-1}\left|\tilde{m}_0\left(\frac{2^{-j}}{2\pi}\right)\right|^2
\sup_{\eta_1 \in [0,1/2\pi]}\prod_{j=0}^{\infty}\left|\tilde{m}_0\left(2^{-j-J_1}\eta_1\right)\right|^2,
\]
Since
\[
|\tilde{m}_0(\xi_1)|^2 \le 1+\left(2\pi\sum_{n} |h(n)||n|\right)|\xi_1| \le e^{\left(2\pi\sum_{n} |h(n)||n|\right)|\xi_1|}.
\]
we can continue by concluding that
\[
|\hat{\phi}(\xi_1)|^2 \le C_2 \cdot |2\pi\xi_1|^{-2\gamma}\cdot \prod_{j=0}^{J_1-1}\left|\tilde{m}_0\left(\frac{2^{-j}}{2\pi}\right)\right|^2
\cdot e^{2^{-J_1+1}\sum_{n} |h(n)||n|}.
\]
The proposition is proved.
\end{proof}

We can conclude from this result that $|\hat{\phi}(\xi_1)| = {O}(|\xi_1|^{-\gamma})$ for $0 \le K' \le \frac{4K-L+1}{5}$
(which is equivalent to $\frac{K-K'}{K'+L-1} \ge \frac14$), where $\gamma = (K-K')-\frac{1}{2}\log_2(C_2)$. Let us now take
a closer look at the decay rate $\gamma$. This depends on $K'$ for fixed $K$ and $L$, where $K'$ can be chosen to be an
arbitrary integer satisfying $0 \leq K' < K$. Now regarding the parameters $K'$ and $K$ as control parameters for the low pass filter $m_0$ given by
\[
|m_0(\xi_1)|^2 = {\cos(\pi\xi_1)}^{2(K-K')} \left(   {\cos(\pi\xi_1)}^{2(K')}\sum_{n=0}^{L-1}{K-1+n\choose n}(\sin(\pi\xi_1))^{2n} \right),
\]
we see that, for any choice of $K'$ satisfying $0 \leq K' < K$, the same scaling function $\phi$ is generated. Therefore,
for analyzing the decay of $\hat{\phi}$, it is sufficient to show that there exists some $K'$ such that $\gamma$ is above
a certain threshold; later we aim for $\gamma > 3$. The following lemma makes this consideration explicit by choosing $K' = 0$.
We however would like to emphasize that this choice was made only for proving feasibility, i.e., \eqref{eq:psi}, of the shearlet which we
will introduce in Theorem \ref{theo:compact}; it will not be our optimal choice when estimating $L_{sup}$ and $R(c)$
of the generated shearlet frame. The technical proof of this lemma is deferred to Subsection \ref{subsec:compactproof}.

\begin{lemma}\label{lemma:gamma}
Let $m_0$ be the low pass filter defined in \eqref{eq:m0} with $K \ge \frac{3L}{2}$ and $L \ge 2$,
and let $C_2$  and $\tilde{m}_0$ be as in Lemma \ref{lemma:tildem} with $K'=0$. Then
\[
\max_{\xi_1 \in [0,1]}|\tilde{m}_0(\xi_1)|^2 \le 2^{2K-L/2-1}
\]
and the constant $\gamma$ from Proposition \ref{proposition:piUpper}
satisfies.
\[
\gamma = K-\frac12\log_2(C_2) > \frac12(\frac L2 +1).
\]
\end{lemma}

\begin{figure}[ht]
\begin{center}
\includegraphics[height=1.2in]{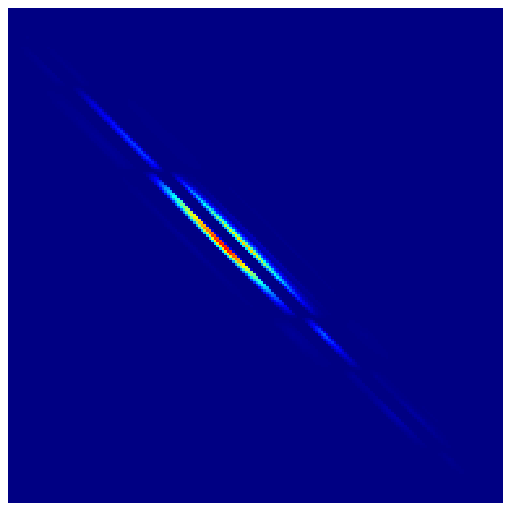}
\hspace*{-1cm}
\includegraphics[height=1.2in]{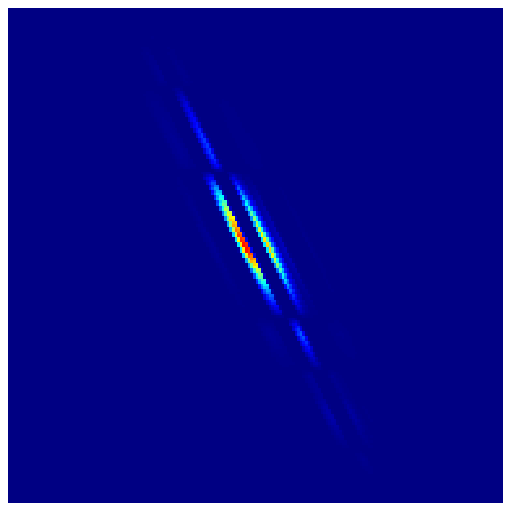}
\hspace*{-1cm}
\includegraphics[height=1.2in]{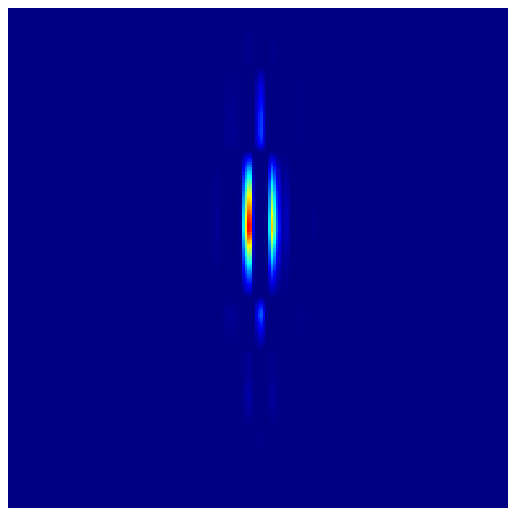}
\hspace*{-1cm}
\includegraphics[height=1.2in]{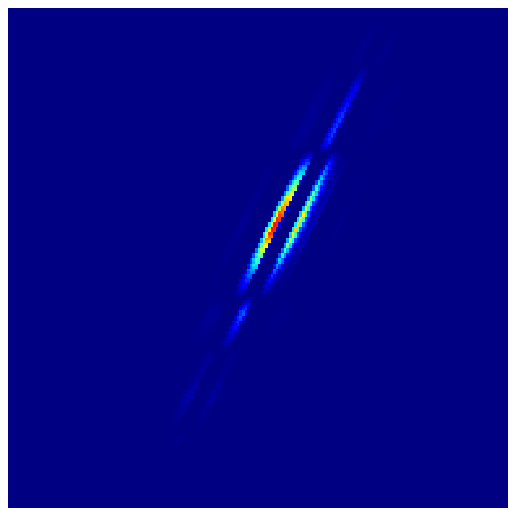}
\hspace*{-1cm}
\includegraphics[height=1.2in]{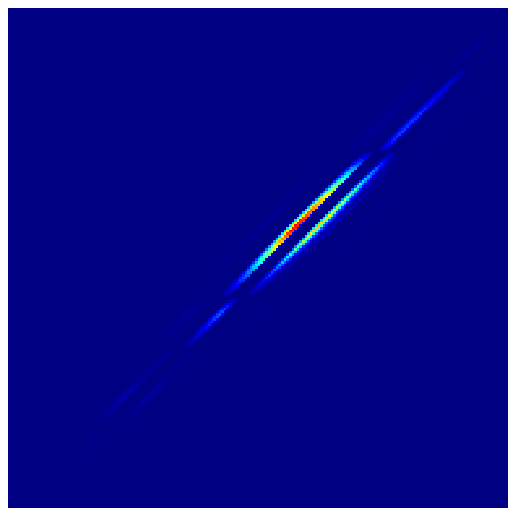}
\end{center}
\caption{Display of compactly supported shearlets $\psi_{2,k,0}, k = -2,-1,0,1,2$, in spatial domain. }
\label{fig:compsuppshearlets}
\end{figure}

Finally, we have reached the stage, where we can introduce a compactly supported shearlet $\psi \in L^2(\RR^2)$,
which generates a  (separable) feasible shearlet frame $\mathcal{S}\mathcal{H}(\Lambda,\psi)$ with
`good' mathematical controllable frame bounds. For
illustrative purposes, some
elements of the shearlet frame $\mathcal{S}\mathcal{H}(\Lambda,\psi)$ are displayed in Figure \ref{fig:compsuppshearlets}.

\begin{proposition}\label{proposition:upperlowershear}
Let $K, L \in \bZ^+$ be such that $L \ge 10$ and $ \frac{3L}{2} \le K \le 3L-2$, let $m_0$ be the associated
low pass filter as defined in \eqref{eq:m0}, and let $\phi$ be the associated scaling function as
defined in \eqref{eq:phi}. Further, we define the bandpass filter $m_1$ by
\[
|m_1(\xi_1)|^2 = |m_0(\xi_1+1/2)|^2, \quad \xi_1 \in \RR.
\]
Then the shearlet defined by
\[
\hat{\psi}(\xi) = m_1(4\xi_1)\hat{\phi}(\xi_1)\hat{\phi}(2\xi_2), \quad \xi = (\xi_1,\xi_2) \in \RR^2
\]
satisfies the feasibility condition \eqref{eq:psi}; in particular, for any $0 \le K' \le \frac{4K-L+1}{5}$, we have
\begin{equation}\label{eq:uppershearlet}
|\hat{\psi}(\xi)| \le \min\{1,|q\xi_1|^{\alpha}\}\min\{1,|q'\xi_1|^{-\gamma}\}\min\{1,|r\xi_2|^{-\gamma}\},
\end{equation}
where
\begin{equation}\label{eq:const2}
\alpha = K-K'\quad \mbox{and} \quad \gamma = (K-K')-\frac12\log_2(C_2)
\end{equation}
with $C_2$ being defined in \eqref{eq:c2}, and
\begin{equation}\label{eq:const1}
q = 4\pi C_2^{\frac{1}{2(K-K')}}, q' = 2\pi \left(C_2\prod_{j=0}^{J_1-1}
\left|\tilde{m}_0\left(\frac{2^{-j}}{2\pi}\right)\right|^2
e^{2^{-J_1+1}\sum_{n} |h(n)||n|}\right)^{-\frac{1}{2\gamma}} \hspace*{-0.5cm}, \,\,r=2q'.
\end{equation}
Further, we have
\begin{equation}\label{eq:lowershearlet}
|\hat{\psi}(\xi)|^2
\ge |m_0(\tfrac16)|^2 \cdot \left(\prod_{j=0}^{J_0-1}\left|m_0\left(\frac{2^{-j}}{6}\right)\right|^2e^{-2^{-J_0+2}(1-|m_0(\frac16)|^2)}\right)^2
\cdot \chi_{\Omega}(\xi)> 0, \quad \xi \in \RR^2,
\end{equation}
where $\Omega = \{\xi=(\xi_1,\xi_2) \in \bR^2 \,:\, \xi_1 \in [\frac{1}{12},\frac{1}{6}]\cup[-\frac{1}{12},-\frac{1}{6}],\,\, \xi_2 \in [-\frac{1}{12},\frac{1}{12}]\}$.
\end{proposition}

Now for $K, L \in \bZ^+$ satisfying $L \ge 10$ and $\frac{3L}{2} \le K \le 3L-2$, Proposition \ref{proposition:upperlowershear}
together with \eqref{eq:Upp2} implies that $L_{sup} < \infty$. Moreover, by Proposition
\ref{prop:estimateR}, $R(c)$ can be made arbitrarily small by choosing a sampling matrix $M_c$
with sufficiently small determinant. Finally, let the set $\Omega$ be defined as in Proposition
\ref{proposition:upperlowershear}.
From \eqref{eq:lowershearlet}, we see that $|\hat{\psi}(\xi)|^2 \ge \tilde{L}_{inf} \cdot \chi_{\Omega}(\xi)$, where
\begin{equation}\label{eq:lower_lower}
\tilde{L}_{inf} := \left|m_0\left(\tfrac16\right)\right|^2\left(\prod_{j=0}^{J_0-1}
\left|m_0\left(\frac{2^{-j}}{6}\right)\right|^2e^{-2^{-J_0+2}(1-|m_0(\frac16)|^2)}\right)^2
\end{equation}
Since
\[
\bigcup_{j=0}^{\infty}\bigcup_{|k| \le \lceil 2^{j/2} \rceil} A_{2^j}S^T_{k}\Omega = \cC,
\]
it then follows that
\[
\Phi(\xi,0) = \sum_{j,k}|\hat{\psi}(S_k^{T}A_{2^{-j}\xi})|^2 \ge \tilde{L}_{inf} \cdot
\sum_{j,k}\chi_{\Omega}(S_k^{T}A_{2^{-j}}\xi) \ge \tilde{L}_{inf} \quad \text{on}\quad\cC.
\]
Thus, by Theorem \ref{theo:general}, both the lower and upper frame bound exist with explicit estimates
for the shearlet frame $\mathcal{S}\mathcal{H}(\Lambda,\psi)$ for a sampling matrix $M_c$
with sufficiently small determinant. This directly implies the following main result.

\begin{Thm} \label{theo:compact}
Let
\[
\Lambda = \{(2^{-j},k2^{j/2},S_{k2^{j/2}}A_{2^{j}}M_cm) : j \ge 0,|k| \le \lceil 2^{j/2} \rceil, m \in \bZ^2\}
\]
be the regular sampling set of $\bS_{cone}$.
Further, let $K, L \in \bZ^+$ be such that $L \ge 10$ and $\frac{3L}{2} \le K \le 3L-2$, and define a shearlet $\psi \in L^2(\RR^2)$ by
\[
\hat{\psi}(\xi) = m_1(4\xi_1)\hat{\phi}(\xi_1)\hat{\phi}(2\xi_2), \quad \xi = (\xi_1,\xi_2) \in \RR^2,
\]
where $m_0$ is the low pass filter satisfying
\[
|m_0(\xi_1)|^2 = (\cos(\pi\xi_1))^{2K}\sum_{n=0}^{L-1} {K-1+n\choose n}(\sin(\pi\xi_1))^{2n}, \quad \xi_1 \in \RR,
\]
$m_1$ is the associated bandpass filter defined by
\[
|m_1(\xi_1)|^2 = |m_0(\xi_1+1/2)|^2, \quad \xi_1 \in \RR,
\]
and $\phi$ is the scaling function given by
\[
\hat{\phi}(\xi_1) = \prod_{j=0}^{\infty} m_0(2^{-j}\xi_1), \quad \xi_1 \in \RR.
\]
Then there exists a sampling constant $\hat c_1>0$ such that the shearlet system $\mathcal{S}\mathcal{H}(\Lambda,\psi)$
forms a frame for $L^2(\cC)$ for any sampling matrix $M_c$ with $c=(c_1,c_2) \in (\RR^+)^2$ and $c_2 \le c_1 \le \hat{c}_1$.
Furthermore, the corresponding frame bounds $A$ and $B$ satisfy
\[
\frac{1}{|\det(M_c)|} [\tilde{L}_{inf} - R(c)] \le A \le B \le \frac{1}{|\det(M_c)|} [\tilde{L}_{sup} + R(c)],
\]
where $R(c) < \tilde{L}_{inf} \le L_{inf}$ and $L_{sup} \le \tilde{L}_{sup}$.
\end{Thm}

To illustrate how this result can be applied, we discuss the derived estimates for the frame bounds
using the following particular choices for the free parameters $K, K', L, c_1$, and $c_2$. However,
we do not claim that this is the optimal choice. In fact, optimizing the estimate for the ratio
of the frame bounds $A/B$ is a highly delicate task in this situation.

\begin{example}
For our discussion, we decided upon the following choices for the various parameters necessary for
exploiting Theorem \ref{theo:compact} to derive an estimate for the ratio of the frame bounds:

We may choose $\tilde{L}_{inf}$ as in \eqref{eq:lower_lower} for the expression for the lower frame bound.
The value $\tilde{L}_{sup}$ can be estimated using \eqref{eq:Upp2} with $\mu =\frac{1}{2}, p=1$.
Further, for a given sampling matrix $M_c$ with $c = (c_1,c_2) \in (\RR^+)^2$ and $c_1 \ge c_2 > 0$, the value
$R(c)$ can be estimated using \eqref{eq:upper_R} with $\mu =\frac 12, p = 1$ from Proposition \ref{prop:estimateR}.
For both those estimates, given positive integers $K, K',$ and $L$ such that
for $K' \le \frac{4K-L+1}{5}$ we choose $\mu = \tfrac12$, the values $q, q', r$ are chosen as in \eqref{eq:const1},
and the values $\alpha, \gamma$ as in  \eqref{eq:const2}.

Table \ref{tab:cs}
shows some estimates for the ratio of the frame bounds of $\mathcal{S}\mathcal{H}(\Lambda,\psi)$
using Theorem \ref{theo:compact} for specific choices of the parameters $K, K',$ and $L$ with various sampling constants $c_1, c_2$.
It should be pointed out that numerical experiments show much better results, but the presented ones are the estimates we are
able to {\em prove theoretically}. Such estimates do not even exist for compactly supported wavelet frames!

\begin{table}\label{tab:cs}
\begin{tabular}{|c|c|c|c|c|c|}
\hline
$K$ & $L$ & $c_1$ & $c_2$ & $K'$ & $B/A$   \\ \hline
39 & 18 & 1.00 & 0.40  & (27,17) & 37.1204 \\
39 & 18 & 1.00 & 0.30  & (27,15) & 32.0208 \\
39 & 18 & 1.00 & 0.25 & (27,15) & 31.9105 \\
39 & 18 & 1.00 & 0.15 & (27,15) & 31.9019 \\
39 & 18 & 1.00 & 0.10  & (27,15) & 31.9019 \\
\hline
\end{tabular}
\begin{tabular}{|c|c|c|c|c|c|}
\hline
$K$ & $L$ & $c_1$ & $c_2$ & $K'$ & $B/A$ \\ \hline
39 & 19 & 0.9 & 0.40  & (27,18) & 44.5359 \\
39 & 19 & 0.9 & 0.30  & (27,16) & 28.4307 \\
39 & 19 & 0.9 & 0.25 & (27,15) & 28.0983 \\
39 & 19 & 0.9 & 0.20  & (27,15) & 28.0699 \\
39 & 19 & 0.9 & 0.15 & (27,15) & 28.0683 \\
\hline
\end{tabular}\\[1ex]
\hspace*{0.5cm}{\rm (a)}\hspace*{7.5cm}{\rm (b)}\\
\vspace*{0.3cm} \caption{Some estimates for the ratio of the frame bounds for the compactly
supported shearlet frame constructed in Theorem \ref{theo:compact} for various choices of
the parameters  $K, K',L$ and the sampling constants $c_1, c_2$.
$K'$ is given in the form $(\cdot,\cdot)$, where the first component is an estimate for $L_{sup}$
chosen by \eqref{eq:Upp2} and the second component is an estimate for $R(c)$
chosen by \eqref{eq:upper_R}.}
\end{table}
Notice the flexibility in choosing $K, K',$ and $L$. A full blown optimization analysis might lead to
an even better estimate for the ratio of the frame bounds of $\mathcal{S}\mathcal{H}(\Lambda,\psi)$,
but this is beyond the scope of this paper.
\end{example}

\subsection{Shearlet Frame for $L^2({\bR}^2)$}
\label{subsec:generalRR}

We now aim to generate a frame for the whole space $L^2(\RR^2)$ using the compactly supported (separable)
shearlet frame introduced in the previous subsection. One classical way for achieving this is by
orthogonally projecting each shearlet element in Fourier domain onto the cone $\cC$, and proceeding
likewise for the vertical cone. However, this procedure is quite counterproductive in the sense that
it destroys most of the advantageous properties, e.g., compact support and regularity, of the
elements of such a shearlet frame, in particular, if the Fourier transform of the to be analyzed signal
is not entirely supported in either the horizontal or the vertical cone.

This problem can though be quite easily resolved by taking the union of two cone-adapted shearlet
systems -- one for the horizontal cone $\cC = \cC_1 \cup \cC_3$ and the other for the vertical cone
$\widetilde{\cC} = \cC_2 \cup \cC_4$. In fact, a shearlet frame $\cSH(c;\phi,\psi,\tilde{\psi})$
(cf. Definition \ref{defi:discreteshearlets}) for $L^2(\RR^2)$ can be constructed by the following

\begin{theorem} \label{theo:completeframe}
Let $\psi \in L^2(\RR^2)$ be the shearlet with associated scaling function $\phi \in L^2(\RR)$ both
introduced in Theorem \ref{theo:compact}, and set $\phi(x_1,x_2)=\phi(x_1)\phi(x_2)$ and
$\tilde{\psi}(x_1,x_2) = \psi(x_2,x_1)$.
Then the corresponding shearlet system $\cSH(c;\phi,\psi,\tilde{\psi})$ forms a frame for $L^2(\bR^2)$
for any sampling matrices $M_c$ and $\tilde{M}_c$ with $c=(c_1,c_2) \in (\RR^+)^2$ and $c_2 \le c_1 \le \hat{c}_1$.
\end{theorem}

\begin{proof}
In the sequel, the constants $\alpha, \gamma$ and $q,q',r$ will be those defined in \eqref{eq:const2} and \eqref{eq:const1},
respectively.

First note that the function $\Phi$ defined in \eqref{eq:Phi} now becomes $\Phi : \bR^2 \times \bR^2 \to \bR$ defined by
\begin{equation}\label{eq:moPhi}
\Phi(\xi,\omega) = |\hat \theta(\xi)||\hat \theta(\xi+\omega)|+ \Phi_1(\xi,\omega)+\Phi_2(\xi,\omega),
\end{equation}
where
\[
\Phi_1(\xi,\omega) = \sum\limits_{j \ge 0} \sum\limits_{|k| \le \lceil 2^{j/2} \rceil}
\left|\hat{\psi}(S_{{k}}^{T} A_{2^{-j}}\xi) \right| \left|\hat{\psi}({S_{k}}^{T}A_{2^{-j}}\xi + \omega) \right|
\]
and
\[
\Phi_2(\xi,\omega)=\sum\limits_{j \ge 0} \sum\limits_{|k| \le \lceil 2^{j/2} \rceil}
\left|\hat{\tilde{\psi}}({S}_{{k}} \tilde{A}_{2^{-j}}\xi) \right|
\left|\hat{\tilde{\psi}}({S}_{{k}}\tilde{A}_{2^{-j}}\xi+ \omega) \right|.
\]
Also, $R(c)$ defined in  \eqref{eq:defR}, is now given by the term
\[
\sum_{m \neq 0} \left(\Gamma_0(c_1^{-1}m)\Gamma_0(-c_1^{-1}m)\right)^{\frac 12}
+\left(\Gamma_1(M_c^{-1}m)\Gamma_1(-M_c^{-1}m)\right)^{\frac 12}
+(\Gamma_2(\tilde{M}_c^{-1}m)\Gamma_2(-\tilde{M}_c^{-1}m))^{\frac 12},
\]
where
\[
\Gamma_0(\omega) = \esssup_{\xi \in \RR^2} |\hat \theta(\xi)||\hat \theta(\xi+\omega)| \quad \text{and} \quad
\Gamma_{i}(\omega) = \esssup_{\xi \in \RR^2} \Phi_{i}(\xi,\omega) \quad \text{for}\,\, i = 1,2.
\]
In fact, Theorem \ref{theo:general} can be easily extended to this case so that we only need to estimate
$L_{inf}$ by some positive $\tilde{L}_{inf}$, $L_{sup}$ by some finite $\tilde{L}_{sup}$, and $R(c)$ to derive finiteness of the
estimate for the ratio of the frame bounds in \eqref{eq:estimateAB}. Using the
function $\Phi(\xi,\omega)$ from \eqref{eq:moPhi}, $L_{inf}$ and $L_{sup}$ are defined as in \eqref{eq:defiLinLsup} and \eqref{eq:defR}
-- note that in this case, $\cC$ is replaced by $\RR^2$ in the definition.

We start by estimating the first term $|\hat \theta(\xi)||\hat \theta(\xi+\omega)|$ from \eqref{eq:moPhi}.
WLOG, we may assume that $\|m\|_{\infty} = m_1 \ne 0$ for $m \in \bZ^2 \backslash \{0\}$. Then, by Proposition \ref{proposition:piUpper},
we have
\begin{eqnarray}\label{eq:theta}
|\hat{\theta}(\xi)||\hat{\theta}(\xi+c_1^{-1}m)| &\le& |\hat{\phi}(\xi_1)||\hat{\phi}(\xi_1+c_1^{-1}m)| \nonumber \\
&\le& \min\{1,|q'\xi_1|^{-\gamma}\}\min\{1,|q'\xi_1+q'c_1^{-1}m_1|^{-\gamma}\} \nonumber \\
&\le& \left( \frac{2c_1}{q'}\right)^{\gamma}\|m\|_{\infty}^{-\gamma},
\end{eqnarray}
provided that $c_1^{-1}q'm_1 \ge 4$.
Now choose $c_1>0$ sufficiently small so that $c_1^{-1}q' \ge 4.$

For the second and third term in \eqref{eq:moPhi}, we can apply the estimates \eqref{eq:Upp2} and \eqref{eq:decayR} from the proof of
Proposition \ref{prop:estimateR}. In particular for the third term, we only need to switch the roles of the variables $\xi_1$ and $\xi_2$
to reach the same estimates as \eqref{eq:Upp2} and \eqref{eq:decayR}.

Set now $\gamma'' = \gamma-\gamma'$ for an arbitrarily fixed $\gamma'$ satisfying $1 < \gamma' < \gamma-2$. Then, by \eqref{eq:theta} and \eqref{eq:decayR},
we obtain
\begin{equation}\label{eq:moR}
R(c) \le C\left( \frac{2c_1}{q'}\right)^{\gamma}+C'\left(\frac{2qc_2}{q'r}\right)^{\gamma''}.
\end{equation}
Also, using \eqref{eq:Upp2},
\begin{equation}\label{eq:moUpp}
L_{sup} \le 1+2 \Bigl(  \frac{q}{r}\left( 2+\frac{2}{2\gamma-1}\right)+1 \Bigr)
\Bigl(\Big\lceil \log_{2}\Big(\frac{q}{q'}\Big) \Big\rceil +\frac{1}{1-2^{1-2\alpha}}+\frac{1}{1-2^{-2\gamma}}\Bigr).
\end{equation}
since $|\hat{\theta}(\xi)| \le 1$. This gives the upper bound of $L_{sup}$ and we may choose the upper bound $\tilde{L}_{sup}$
to be the right hand side of \eqref{eq:moUpp}.

Finally, let us estimate the lower bound of $L_{inf}$. By \eqref{eq:lowerphi} and \eqref{eq:lowershearlet}, we have
\beq \label{eq:lowerbounds}
|\hat \theta(\xi)|^2 > C \chi_{\Omega_0}(\xi),\quad |\hat \psi(\xi)|^2>C \chi_{\Omega_1}(\xi), \quad \text{and} \quad
|\hat{\tilde{\psi}}(\xi)|^2 > C \chi_{\Omega_2}(\xi),
\eeq
where
\begin{eqnarray*}
\Omega_0 &=& \{\xi \in \bR2 \,:\, \|\xi\|_{\infty} \le \tfrac16\}, \\
\Omega_1 &=& \{\xi \in \bR^2 \,:\, \tfrac{1}{12} < |\xi_1| < \tfrac16, |\xi_2| < \tfrac{1}{12} \}, \\
\Omega_2 &=& \{\xi \in \bR^2 \,:\, \tfrac{1}{12} < |\xi_2| < \tfrac16, |\xi_1| < \tfrac{1}{12} \}.
\end{eqnarray*}
Setting
\[
\Omega = \left( \Omega_0 \cup \bigcup_{j,k} A_{2^j}S_{k} \Omega_1 \cup \bigcup_{j,k} \tilde{A}_{2^j}S^T_{k} \Omega_2\right),
\]
we observe that $\Omega = \bR^2$ and, by \eqref{eq:lowerbounds},
\begin{equation}\label{eq:moLow}
\Phi(\xi,0) > \tilde{L}_{inf} \cdot \chi_{\Omega}(\xi) \quad \mbox{for some} \,\, \tilde{L}_{inf}>0.
\end{equation}
This implies $\Phi(\xi,0) > \tilde{L}_{inf} >0$ on $\bR^2$, which yields the lower bound of $L_{inf}$.

Concluding, the estimates \eqref{eq:moR}, \eqref{eq:moUpp}, and \eqref{eq:moLow} provide all required constants
$\tilde{L}_{inf}$, $\tilde{L}_{sup}$ and the upper bound of $R(c)$ in \eqref{eq:estimateAB}.
Especially, we can choose a sampling matrix $M_c$ with sufficiently small determinant such that $\tilde{L}_{inf}-R(c)>0$.
This proves the existence of the frame bounds.

Finally, it is obvious that all functions in $\mathcal{S}\mathcal{H}(\Lambda,\tilde{\Lambda},\psi,\tilde \psi, \theta)$ are compactly supported in
spatial domain, which finished the proof.
\end{proof}
In the proof of Theorem \ref{theo:completeframe}, we see that our upper bounds of $L_{sup}$ and $R(c)$ are about twice as
large as the upper bounds in the cone case and the lower bound $\tilde{L}_{inf}$ is about the same compared to the cone case.
Also $c$ can be chosen so that $R(c)$ is sufficiently small.
This indicates that in this case, our estimate for the ratio of the frame bounds is about twice as large as the estimate for
the ratio $B/A$ in the cone case for sufficiently small determinant of the sampling matrix $M_c$.

\subsection{Sparse Approximation using Compactly Supported Shearlets}
\label{subsec:sparse}

One essential performance criterion for a frame composed of anisotropic elements -- besides the ratio of the frame bounds --
is the approximation rate of curvilinear objects. This viewpoint arose due to the fact that edges are the most prominent
features in images, hence representation systems should in particular provide sparse representations for those.

To stand on solid ground, we first briefly recall the mathematical model of a cartoon-like image introduced
in \cite{CD04}. For $\nu > 0$, the set $STAR^2(\nu)$ is defined to be the set of all $B \subset [0,1]^2$ such
that $B$ is a translate of a set $\{x \in \RR^2 : |x| \le \rho(\theta), x = (|x|,\theta) \mbox{ in polar coordinates}\}$
which satisfies $\sup|\rho^{''}(\theta)| \leq \nu$, $\rho \leq \rho_0 < 1$. Then, $\cE^2(\nu)$
denotes the set of functions $f \in L^2(\RR^2)$ of the form
\[
f = f_0 + f_1 \chi_{B},
\]
where $f_0,f_1 \in C_0^2([0,1]^2)$ and $B \in STAR^2(\nu)$. In \cite{Don01}, it was proven that the optimal approximation rate
for such cartoon-like image models which can be achieved under some restrictions on the representation system as well as on
the selection procedure of the approximating coefficients is
\[
\norm{f-f_N}_2^2 \le C \cdot N^{-2} \quad \mbox{as } N \to \infty,
\]
where $f_N$ is the $N$-term approximation generated by the $N$ largest coefficients in magnitude.

In \cite{KL10}, two of the authors proved that a large class of shearlet frames -- including, in particular, a
significant set of compactly supported shearlet frames -- provide optimally sparse approximations of such
cartoon-like images. For the convenience of the reader, we state the result below. Notice that in
\cite{KL10}, the result was proven for a isotropic sampling matrix, but the extension to our anisotropic
sampling is immediate.

\begin{theorem}[\cite{KL10}] \label{theo:sparse}
Let $c > 0$, and let $\phi, \psi, \tilde{\psi} \in L^2(\RR^2)$ be compactly supported. Suppose
that, in addition, for all $\xi = (\xi_1,\xi_2) \in \RR^2$, the shearlet $\psi$ satisfies
\bitem
\item[(i)] $|\hat\psi(\xi)| \le C \cdot \min\{1,|\xi_1|^{\alpha}\} \cdot \min\{1,|\xi_1|^{-\gamma}\} \cdot \min\{1,|\xi_2|^{-\gamma}\}$ and
\item[(ii)] $\left|\frac{\partial}{\partial \xi_2}\hat \psi(\xi)\right|
\le  |h(\xi_1)| \cdot \left(1+\frac{|\xi_2|}{|\xi_1|}\right)^{-\gamma}$,
\eitem
where $\alpha > 5$, $\gamma \ge 4$, $h \in L^1(\bR)$, and $C$ is a constant, and suppose that
the shearlet $\tilde{\psi}$ satisfies (i) and (ii) with the roles of $\xi_1$ and $\xi_2$ reversed.
Further, suppose that $\cSH(c;\phi,\psi,\tilde{\psi})$ forms a
 frame for $L^2(\RR^2)$.

Then, for any $\nu > 0$, the shearlet frame $\cSH(c;\phi,\psi,\tilde{\psi})$ provides (almost)
optimally sparse approximations of functions $f \in \cE^2(\nu)$, i.e., there exists some $C > 0$ such
that
\[
\|f-f_N\|_2^2 \leq C\cdot {(\log{N})}^3\cdot N^{-2}  \qquad \text{as } N \rightarrow \infty,
\]
where $f_N$ is the nonlinear N-term approximation obtained by choosing the N largest shearlet coefficients
of $f$.
\end{theorem}

Using this result, we can prove that in fact the compactly supported shearlet frames we introduced in
Theorem \ref{theo:compact} even provide (almost) optimally sparse approximations of cartoon-like images, i.e.,
functions in $\cE^2$. Although the optimal rate is not achieved, the $\log$-factor is typically considered negligible
compared to the $N^{-2}$-factor, wherefore the term `almost optimal' has been adopted into the
language.

\begin{theorem}
Let $\psi \in L^2(\RR^2)$ be the shearlet with associated scaling function $\phi \in L^2(\RR)$ both
introduced in Theorem \ref{theo:compact} with $L>18$, and set $\phi(x_1,x_2)=\phi(x_1)\phi(x_2)$ as well as
$\tilde{\psi}(x_1,x_2) = \psi(x_2,x_1)$.  Then the compactly supported shearlet frame
$\cSH(c,\tilde{c};\phi,\psi,\tilde{\psi})$
provides (almost) optimally sparse approximations of functions $f \in \cE^2(\nu)$, i.e., there exists
some $C > 0$ such that
\[
\|f-f_N\|_2^2 \leq C\cdot {(\log{N})}^3\cdot N^{-2}  \qquad \text{as } N \rightarrow \infty,
\]
where $f_N$ is the nonlinear N-term approximation obtained by choosing the N largest shearlet coefficients
of $f$.
\end{theorem}

\begin{proof}
The fact that $\cSH(c;\phi,\psi,\tilde{\psi})$ is a compactly supported frame follows from Theorem
\ref{theo:completeframe}. Hence it remains to prove that the shearlet $\psi$ satisfies
conditions (i) and (ii) from Theorem \ref{theo:sparse}, and the shearlet $\tilde{\psi}$ likewise.

By Lemma \ref{lemma:gamma} and Proposition \ref{proposition:upperlowershear}, there exists a
constant $C$ such that
\[
|\hat \psi(\xi)| \le C \cdot \min\{1,|\xi_1|^K\}\min\{1,|\xi_1|^{-L/4-1/2}\}\min\{1,|\xi_2|^{-L/4-1/2}\},
\]
where $L/4+1/2 > 4$ and $K>5$. This already proves condition (i) of Theorem \ref{theo:sparse}.

To show condition (ii), we first observe that there exists a function $\psi_1$ such that
we can write $\hat \psi (\xi) = \hat \psi_1(\xi_1)\hat \phi(2\xi_2)$
so that
\[
|\hat \psi_1(\xi_1)| \leq C\cdot\min\{1,|\xi_1|^K\}\min\{1,|\xi_1|^{-L/4-1/2}\} \quad \text{and}
\quad |\hat\phi(\xi_2)| \leq C\cdot\min\{1,|\xi_2|^{-L/4-1/2}\}.
\]
Hence it is sufficient to prove that
\beq \label{eq:enough}
|(\hat \phi)^{'}(\xi_2)| \leq C \cdot \min\{1,|\xi_2|^{-\gamma'}\}\quad \mbox{with }\gamma'\ge 4.
\eeq
For this, we first choose  a positive integer $\rho'>0$ such that $4 \leq \rho' < \frac{L}{4}-\frac12 < K$, which
we are allowed to do since $L > 18$. Then, by the definition of $\phi$,
\begin{eqnarray*}
\hat \phi(\xi_2) = \left(\frac{\sin (2\pi \xi_2)}{2\pi\xi_2}\right)^{\rho'}
\left(\frac{\sin (2\pi \xi_2)}{2\pi\xi_2}\right)^{K-\rho'}\prod_{j=0}^{\infty}\tilde{m}_0(2^{-j}\xi_2),
\end{eqnarray*}
where $|\tilde{m}_0(\xi_2)|^2 = \sum_{n=0}^{L-1}{K-1+n\choose n}(\sin(\pi \xi_2))^{2n}.$
Now define $S$ and $\tilde{\phi}$ by
\[
S(\xi_2) = \left(\frac{\sin (2\pi \xi_2)}{2\pi\xi_2}\right)^{\rho'}
\quad \text{and} \quad
\hat{\tilde{\phi}}(\xi_2)=\left(\frac{\sin (2\pi \xi_2)}{2\pi\xi_2}\right)^{K-\rho'}\prod_{j=0}^{\infty}\tilde{m}_0(2^{-j}\xi_2),
\]
which gives
\[
(\hat{\phi})^{'}(\xi_2) = S'(\xi_2)\hat{\tilde{\phi}}(\xi_2)+S(\xi_2)(\hat{\tilde{\phi}})'(\xi_2).
\]
Since $\max\{|S'(\xi_2)|,|S(\xi_2)|\} \leq C|\xi_2|^{-\rho'}$, it remains to show that both $|\hat{\tilde{\phi}}|$
and $|(\hat{\tilde{\phi}})|$ are bounded. Notice that we may assume that $\tilde{m}_0(0)=1$ and $\tilde{\phi}$ is
compactly supported (see \cite{Dau92}). Therefore, boundedness of $\tilde \phi$ implies that $\hat{\tilde{\phi}}$ and
$(\hat{\tilde{\phi}})'$ are bounded, which enables us to restrict our task further and it remains to check
that $\hat{\tilde{\phi}} \in L^1(\bR)$. This will now be proved by showing that, for all $\xi_2$,
\beq \label{eq:enough2}
|\hat{\tilde{\phi}}(\xi_2)| \leq C \cdot \min\{1,|\xi_2|^{-\eta}\}\quad \mbox{with }\eta >1,
\eeq
which obviously implies $\hat{\tilde{\phi}} \in L^1(\bR).$

First, it is easy to show that $|\hat{\tilde{\phi}}|$ is bounded
on $[-1,1].$ Therefore, if $|\xi_2| \le 1$, then \eqref{eq:enough2} holds.
Now assume that $|\xi_2|>1$. Then there exists a positive integer $J>0$
such that $2^{J-1} \leq |\xi_2| \leq 2^J$, and, using the fact that Lemma \ref{lemma:gamma} implies
$\max_{\xi_2}|\tilde{m}_0(\xi_2)|^2 < 2^{2K-L/2-1}$, we obtain
\begin{eqnarray*}
|\hat{\tilde \phi}(\xi_2)| &\leq& C|\xi_2|^{-(K-\rho')}\prod_{j=0}^{J-1}\max_{\xi_2} |\tilde{m}_0(\xi_2)|
\prod_{j=0}^{\infty}|\tilde{m}_0(2^{-j-J}\xi_2)| \\
&\leq& C|\xi_2|^{-(K-\rho'-K+L/4+1/2)}\sup_{\xi \in [-1,1]}\prod_{j=0}^{\infty}|\tilde{m}_0(2^{-j}\xi_2)|
\\
&\leq& C|\xi_2|^{-(L/4+1/2-\rho')}.
\end{eqnarray*}
Hence, $L/4+1/2-\rho' > 1$ and this proves \eqref{eq:enough2}, and hence \eqref{eq:enough}.
The theorem is proved.
\end{proof}

\section{Proofs}
\label{sec:proofs}

\subsection{Proofs of Results from Section \ref{sec:generalsufficient}}

\subsubsection{Proof of Theorem \ref{theo:general}}
\label{subsubsec:general}

Let $f \in L^{2}(\cC)$. Then, by definition of $\psi_{j,k,m}$,
\begin{eqnarray}\nonumber
\lefteqn{\sum\limits_{j \in \mathbb{Z}} \sum\limits_{k \in K_j} \sum\limits_{m \in \mathbb{Z}^{2}}
\left|\langle \hat{f}, \hat{\psi}_{j,k,m}\rangle_{L^{2}(\cC)} \right|^2}\\ \label{eq:pf_general1}
& = & \sum\limits_{j \ge 0} \sum\limits_{k \in K_j} \sum\limits_{m \in \mathbb{Z}^{2}}
a_{j}^{3/2}\Bigg| \int\limits_{\cC} \hat{f}(\xi) \overline{\hat{\psi}(S_{s_{k}}^{T}A_{a_{j}}\xi )}
e^{2\pi i \langle  \xi, A_{a_{j}}S_{s_{k}} c m \rangle} \, d\xi \, \Bigg|^2.
\end{eqnarray}
We now first decompose the sum over $m \in \bZ^2$. For this, set $\Omega =
\left[-\frac{1}{2}, \frac{1}{2}\right]^{2}$. Then, by appropriate changes of
variables,
{\allowdisplaybreaks
\begin{eqnarray*}
\lefteqn{\sum\limits_{m \in \mathbb{Z}^{2}}
a_{j}^{3/2}\Bigg|\int\limits_{\cC} \hat{f}(\xi) \overline{\hat{\psi}(S_{s_{k}}^{T}A_{a_{j}}\xi )}
e^{2\pi i \langle  \xi, A_{a_{j}}S_{s_{k}} c m \rangle} \, d\xi \, \Bigg|^2}\\[1ex]
& = & \sum\limits_{m \in \mathbb{Z}^{2}} \frac{a_{j}^{-3/2}}{|\det(M_c)|}
\Bigg|\int\limits_{\mathbb{R}^{2}} \hat{f}\left(A_{a_{j}}^{-1}S_{s_{k}}^{-T}M_c^{-1}\xi \right)
\chi_{\cC}\left(A_{a_{j}}^{-1}S_{s_{k}}^{-T}M_c^{-1}\xi \right) \overline{\hat{\psi}\left(M_c^{-1}\xi \right)}
e^{2\pi i \langle \xi , m \rangle} \, d\xi \Big|^2\\[1ex]
& = & \sum\limits_{m \in \mathbb{Z}^{2}} \frac{a_{j}^{-3/2}}{|\det(M_c)|}
\Bigg|\sum_{\ell \in \bZ^2} \int_{\Omega + \ell} \hat{f}\left(A_{a_{j}}^{-1}S_{s_{k}}^{-T}M_c^{-1}\xi \right)
\chi_{\cC}\left(A_{a_{j}}^{-1}S_{s_{k}}^{-T}M_c^{-1}\xi \right) \overline{\hat{\psi}\left(M_c^{-1}\xi \right)}
e^{2\pi i \langle \xi , m \rangle} \, d\xi \Big|^2\\[1ex]
& = & \sum\limits_{m \in \mathbb{Z}^{2}} \frac{a_{j}^{-3/2}}{|\det(M_c)|}
\Bigg|\int_{\Omega} \sum_{\ell \in \bZ^2} \hat{f}\left(A_{a_{j}}^{-1}S_{s_{k}}^{-T}M_c^{-1}(\xi+\ell) \right)
\chi_{\cC}\left(A_{a_{j}}^{-1}S_{s_{k}}^{-T}M_c^{-1}(\xi+\ell) \right)\\
& & \cdot \overline{\hat{\psi}\left(M_c^{-1}(\xi+\ell) \right)} e^{2\pi i \langle \xi , m \rangle} \, d\xi \Big|^2
\end{eqnarray*}
}
Next, applying Plancherel's theorem,
\begin{eqnarray*}
\lefteqn{\sum\limits_{m \in \mathbb{Z}^{2}}
a_{j}^{3/2}\Bigg|\int\limits_{\cC} \hat{f}(\xi) \overline{\hat{\psi}(S_{s_{k}}^{T}A_{a_{j}}\xi )}
e^{2\pi i \langle  \xi, A_{a_{j}}S_{s_{k}} c m \rangle} \, d\xi \, \Bigg|^2}\\[1ex]
& = & \hspace*{-0.3cm} \frac{a_{j}^{-3/2}}{|\det(M_c)|} \int_\Omega
\Bigg|\sum_{\ell \in \bZ^2} \hat{f}\left(A_{a_{j}}^{-1}S_{s_{k}}^{-T}M_c^{-1}(\xi+\ell) \right)
\chi_{\cC}\left(A_{a_{j}}^{-1}S_{s_{k}}^{-T}M_c^{-1}(\xi+\ell) \right) \overline{\hat{\psi}\left(M_c^{-1}(\xi+\ell) \right)}\Big|^2 d\xi
\end{eqnarray*}
Resolving the absolute values yields,
{\allowdisplaybreaks
\begin{eqnarray*}
\lefteqn{\sum\limits_{m \in \mathbb{Z}^{2}}
a_{j}^{3/2}\Bigg|\int\limits_{\cC} \hat{f}(\xi) \overline{\hat{\psi}(S_{s_{k}}^{T}A_{a_{j}}\xi )}
e^{2\pi i \langle  \xi, A_{a_{j}}S_{s_{k}} c m \rangle} \, d\xi \, \Bigg|^2}\\[1ex]
& = & \frac{a_{j}^{-3/2}}{|\det(M_c)|} \int_\Omega \sum_{m, \ell \in \bZ^2}
\hat{f}\left(A_{a_{j}}^{-1}S_{s_{k}}^{-T}M_c^{-1}(\xi+\ell) \right)
\chi_{\cC}\left(A_{a_{j}}^{-1}S_{s_{k}}^{-T}M_c^{-1}(\xi+\ell) \right)
\overline{\hat{\psi}\left(M_c^{-1}(\xi+\ell) \right)} \\
& & \cdot\overline{\hat{f}\left(A_{a_{j}}^{-1}S_{s_{k}}^{-T}M_c^{-1}(\xi+m) \right)}
\hat{\psi}\left(M_c^{-1}(\xi+m) \right) d\xi\\[1ex]
& = & \frac{a_{j}^{-3/2}}{|\det(M_c)|} \sum_{\ell \in \bZ^2} \int_{\Omega+\ell}
\hat{f}\left(A_{a_{j}}^{-1}S_{s_{k}}^{-T}M_c^{-1}\xi \right)\chi_{\cC}\left(A_{a_{j}}^{-1}S_{s_{k}}^{-T}M_c^{-1}\xi \right)
\overline{\hat{\psi}\left(M_c^{-1}\xi \right)}\\
& & \cdot \sum_{m \in \bZ^2} \overline{\hat{f}\left(A_{a_{j}}^{-1}S_{s_{k}}^{-T}M_c^{-1}(\xi+m-\ell) \right)}
\hat{\psi}\left(M_c^{-1}(\xi+m-\ell) \right) d\xi\\[1ex]
& = & \frac{a_{j}^{-3/2}}{|\det(M_c)|} \int_{\bR^2} \sum_{m \in \bZ^2}
\hat{f}\left(A_{a_{j}}^{-1}S_{s_{k}}^{-T}M_c^{-1}\xi \right)\chi_{\cC}\left(A_{a_{j}}^{-1}S_{s_{k}}^{-T}M_c^{-1}\xi \right)
\overline{\hat{f}\left(A_{a_{j}}^{-1}S_{s_{k}}^{-T}M_c^{-1}(\xi+m) \right)  }\\
& & \cdot  \overline{\hat{\psi}\left(M_c^{-1}\xi \right)}
\hat{\psi}\left(M_c^{-1}(\xi+m) \right) d\xi.
\end{eqnarray*}
}
Combining with \eqref{eq:pf_general1},
{\allowdisplaybreaks
\begin{eqnarray}\nonumber
\sum\limits_{j \in \mathbb{Z}} \sum\limits_{k \in K_j} \sum\limits_{m \in \mathbb{Z}^{2}}
\left|\langle \hat{f}, \hat{\psi}_{j,k,m}\rangle_{L^{2}(C)} \right|^2
& = & \frac{1}{|\det(M_c)|}\int\limits_{\cC} \sum\limits_{j \in \mathbb{Z}} \sum\limits_{k \in K_j}
\sum\limits_{m \in \mathbb{Z}^{2}}\hat{f}\left(\xi \right) \overline{\hat{f}\left(\xi +A_{a_{j}}^{-1}S_{s_{k}}^{-T}M_c^{-1}m\right)}\\ \nonumber
& & \hspace*{2.25cm}\cdot \overline{\hat{\psi} \left(S_{s_{k}}^{T}A_{a_{j}}\xi\right)} \hat{\psi}\left(S_{s_{k}}^{T}A_{a_{j}}\xi + M_c^{-1}m \right)\, d\xi\\[1ex]
& = & T_1 + T_2, \label{eq:pf_general2}
\end{eqnarray}
}
where
\[
T_1 = \frac{1}{|\det(M_c)|}\sum\limits_{j \in \mathbb{Z}} \sum\limits_{k \in K_j} \int\limits_{\cC} \left|\hat{f}\left(\xi \right) \right|^{2}
\left|\hat{\psi} \left(S_{s_{k}}^{T}A_{a_{j}}\xi \right)\right|^{2}\,d\xi
\]
and
\begin{eqnarray*}
T_2 & =  & \frac{1}{|\det(M_c)|}\sum\limits_{j \in \mathbb{Z}} \sum\limits_{k \in K_j}\int\limits_{\cC}\hspace*{-0.1cm}\sum\limits_{m \in \Z^{2}\setminus\{0\}}
\hat{f}\left(\xi \right) \overline{\hat{f}\left(\xi+A_{a_{j}}^{-1}S_{k}^{-T}M_c^{-1}m\right)}\\
& & \hspace*{2.25cm}\cdot \overline{\hat{\psi}
\left(S_{s_{k}}^{T}A_{a_{j}}\xi\right)} \hat{\psi}\left(S_{s_{k}}^{T}A_{a_{j}}\xi+M_c^{-1}m\right)d\xi.
\end{eqnarray*}
By using Cauchy-Schwartz inequality,
\beq \label{eq:pf_general3}
|T_2| \le \frac{1}{|\det(M_c)|} \|\hat{f}\|^{2} \sum\limits_{m \in \Z^{2}\setminus\{0\}}
\left[\Gamma\left(M_c^{-1}m\right)\Gamma\left(-M_c^{-1}m\right)\right]^{1/2}.
\eeq
By \eqref{eq:pf_general2} and \eqref{eq:pf_general3}, we finally obtain
\[
\frac{L_{inf} - R(c)}{|\det(M_c)|}\|\hat{f}\|^{2}
\le \sum\limits_{j \in \mathbb{Z}} \sum\limits_{k \in K_j} \sum\limits_{m \in \mathbb{Z}^{2}}
\left|\langle \hat{f}, \hat{\psi}_{j,k,m}\rangle_{L^{2}(C)} \right|^2
\le \frac{L_{sup} - R(c)}{|\det(M_c)|}\|\hat{f}\|^{2}.
\]
The introduction of the lower bound for $L_{inf}$ in the estimate $\frac{1}{|\det(M_c)|}\|\hat{f}\|^{2} \left[L_{inf} - R(c) \right]$
is immediate. The theorem is proved.


\subsubsection{Proof of Proposition \ref{prop:estimateR}}
\label{subsubsec:estimateR}

We start by estimating $\Gamma(2\omega_{1},2\omega_{2})$, and will use this later to derive the
claimed upper estimate for $R(c)$. For each $(\omega_{1}, \omega_{2}) \in {\mathbb{R}^{2}}\backslash \{0\}$,
we first split the sum over $j$ as
{\allowdisplaybreaks
\begin{eqnarray}\nonumber
\lefteqn{\Gamma(2\omega_{1},2\omega_{2})}\\[1ex] \nonumber
& = & \esssup_{\xi \in {\mathbb{R}^{2}}} \; \sum\limits_{j \ge 0} \sum\limits_{k \in K_j}
\left|\hat{\psi}\left(a_{j}\xi_{1}, s_{k}a_{j}\xi_{1}+a_{j}^{1/2}\xi_{2}\right)\right|
 \left|\hat{\psi}\left(a_{j}\xi_{1} + 2\omega_{1}, s_{k}a_{j}\xi_{1}+a_{j}^{1/2}\xi_{2}+2\omega_{2} \right)\right|\\[1ex] \nonumber
& \le & \esssup\limits_{\xi \in {\mathbb{R}^{2}}} \Bigg(\sum\limits_{\{j: |a_{j}\xi_{1}| <\|\omega\|_{\infty}\}}
+ \sum\limits_{\{j: |a_{j}\xi_{1}| \ge \|\omega\|_{\infty}\}} \Bigg)\sum\limits_{k \in \Z} \left|\hat{\psi}\left(a_{j}\xi_{1},
s_{k}a_{j}\xi_{1}+a_{j}^{1/2}\xi_{2}\right)\right|\\ \nonumber
& & \cdot \left|\hat{\psi}\left(a_{j}\xi_{1} + 2\omega_{1}, s_{k}a_{j}\xi_{1}+a_{j}^{1/2}\xi_{2}+2\omega_{2} \right)\right|\\[1ex] \label{eq:pf_estimateR1}
& = & \esssup\limits_{\xi \in {\mathbb{R}^{2}}} (I_{1}+I_{2}),
\end{eqnarray}
where
\[
I_1 =\sum\limits_{\{j  :  |a_{j}\xi_{1}| \le \|\omega\|_{\infty} \}}
\sum\limits_{k \in \Z} \left|\hat{\psi}\left(a_{j}\xi_{1},s_{k}a_{j}\xi_{1}+a_{j}^{1/2}\xi_{2}\right)\right|
\left|\hat{\psi}\left(a_{j}\xi_{1} + 2\omega_{1},s_{k}a_{j}\xi_{1}+a_{j}^{1/2}\xi_{2}+2\omega_{2} \right)\right|
\]
and
\[
I_2 = \sum\limits_{\{j : |a_{j}\xi_{1}| > \|\omega\|_{\infty}\}}
\sum\limits_{k \in \Z} \left|\hat{\psi}\left(a_{j}\xi_{1},s_{k}a_{j}\xi_{1}+a_{j}^{1/2}\xi_{2}\right)\right|
\left|\hat{\psi}\left(a_{j}\xi_{1} + 2\omega_{1},s_{k}a_{j}\xi_{1}+a_{j}^{1/2}\xi_{2}+2\omega_{2}\right)\right|.
\]


The next step consists in estimating $I_1$ and $I_2$. Before we delve in the estimations, we first introduce some
useful inequalities which will be used later. Recall that $\alpha>\gamma > 3$, and $q,q',r$ are positive constants
satisfying $q > q' >0$,  $q > r >0$. Further, let $\gamma'' = \gamma-\gamma'$ for an arbitrarily fixed $\gamma'$
satisfying $1 < \gamma' < \gamma-2$. Then we have the following inequalities for $x,y,z \in \RR$.
The proofs are all elementary, and hence we will skip them.
\begin{eqnarray} \nonumber
\lefteqn{\min\{1,|qx|^{\alpha}\}\min\{1,|q'x|^{-\gamma}\}\min\{1,|ry|^{-\gamma}\}}\\ \label{eq:ineq1}
&\leq& \min\{1,|qx|^{\alpha-\gamma}\}\min\{1,|q'x|^{-\gamma}\}\min\left\{1,| (qx)^{-1}ry|^{-\gamma}\right\},
\end{eqnarray}
\begin{equation}\label{eq:ineq2}
\min\{1,|x|^{-\gamma}\}\min\left\{1,\left| \frac{1+z}{x+y}\right|^{\gamma}\right\} \le 2^{\gamma''}|y|^{-\gamma''}\min\{1,|x|^{-\gamma'}\}
\max\{1,|1+z|^{\gamma''}\},
\end{equation}
\begin{equation}\label{eq:ineq3}
\min\{1,|qx|^{\alpha-\gamma}\}\min\{1,|q'x|^{-\gamma}\}|x|^{\gamma''} \le (q')^{-\gamma''},
\end{equation}
and
\begin{equation}\label{eq:ineq4}
\min\{1,|qx|^{\alpha-\gamma}\}\min\{1,|q'x|^{-\gamma}\}|x|^{\gamma''} \le (q')^{-\gamma''}
\min\{1,|qx|^{\alpha-\gamma+\gamma''}\}\min\{1,|q'x|^{-\gamma'}\}.
\end{equation}
We start with $I_1$. By applying \eqref{eq:psi} and \eqref{eq:ineq1},
\begin{eqnarray} \nonumber
I_1
& \le &  \sum\limits_{\{j :  |a_{j}\xi_{1}| \le \|\omega\|_{\infty} \}} \min\{|qa_{j}\xi_{1}|^{{\alpha - \gamma}},1\}\min
\{|q'a_{j}\xi_{1}|^{-\gamma},1\} \cdot \min \{\left|q(a_{j}\xi_{1} + 2\omega_{1})\right|^{{\alpha - \gamma}},1\}  \\ \nonumber
& & \qquad \cdot \min\{ \left|q'(a_{j}\xi_{1} + 2\omega_{1})\right|^{-\gamma},1\}\sum\limits_{k \in \mathbb{Z}}
\min\left\{\left|\frac rq \left(s_{k}+\frac{\xi_{2}}{a_j^{1/2}\xi_{1}}\right)\right|^{-{\gamma}},1\right\} \\ \label{eq:pf_estimateR2}
&& \qquad \cdot\min\left\{
\left|\frac rq \left[\left(\frac{2\omega_{2}}{a_{j}\xi_{1}}\right) +\left(s_{k} + \frac{\xi_{2}}{a_j^{1/2} \xi_{1}}\right)\right]\right|^{{-\gamma}}
\left|1+\frac{2\omega_{1}}{a_{j}\xi_{1}}\right|^{\gamma},1\right\}.
\end{eqnarray}
We now distinguish two cases, namely $\|\omega\|_{\infty} = |\omega_{1}| \ge |a_{j}\xi_{1}|$ and
$\|\omega\|_{\infty} = |\omega_{2}| \ge |a_{j}\xi_{1}|$. Notice
that these two cases indeed encompass all possible relations between $\omega$ and $\xi_1$.

\medskip

{\em Case 1}: We assume that $\|\omega\|_{\infty} = |\omega_{1}| \ge |a_{j}\xi_{1}|$, hence
$|a_{j}\xi_{1}+2\omega_{1}| \ge |\omega_{1}|$. This implies $\left|q'(a_{j}\xi_{1} + 2\omega_{1})\right|^{-\gamma} \le
\|q'\omega\|_{\infty}^{-\gamma} $. Thus, continuing \eqref{eq:pf_estimateR2},
\begin{eqnarray*}
I_1
& \le & \sum\limits_{\{j :  |a_{j}\xi_{1}| \le \|\omega\|_{\infty} \}} \min\{|qa_{j}\xi_{1}|^{{\alpha - \gamma}},1\}\min
\{|q'a_{j}\xi_{1}|^{-\gamma},1\}\left|q'(a_{j}\xi_{1} + 2\omega_{1})\right|^{-\gamma} \\
& & \qquad \cdot  \sum\limits_{k \in \mathbb{Z}}
\min\left\{\left|\frac rq \left(s_{k}+\frac{\xi_{2}}{a_j^{1/2}\xi_{1}}\right)\right|^{-{\gamma}},1\right\} \\
&\le&
\|q'\omega\|_{\infty}^{-\gamma}\sum\limits_{\{j :  |a_{j}\xi_{1}| \le \|\omega\|_{\infty} \}} \min\{|qa_{j}\xi_{1}|^{{\alpha - \gamma}},1\}\min
\{|q'a_{j}\xi_{1}|^{-\gamma},1\}
\\
& & \qquad \cdot
\frac qr \sum\limits_{k \in \mathbb{Z}}
\frac rq \min\left\{\left|\frac rq \left(s_{k}+\frac{\xi_{2}}{a_j^{1/2}\xi_{1}}\right)\right|^{-{\gamma}},1\right\}.
\end{eqnarray*}
By \eqref{eq:sk2} and Lemma \ref{lemma:estimate_aj}, finally conclude
\begin{eqnarray} \nonumber
I_1
& \le & \frac{\left( \frac qr C(\gamma) \right)}{\|q'\omega\|_{\infty}^{\gamma}}\sum_{j \in \bZ}
\min\{|qa_{j}\xi_{1}|^{{\alpha - \gamma}},1\}\min \{|q'a_{j}\xi_{1}|^{-\gamma},1\}\\ \label{eq:case1}
& \le & p\frac{\left( \frac qr C(\gamma) \right)}{\|q'\omega\|_{\infty}^{\gamma}}
\left( \left\lceil \log_{1/\mu}\left( \frac {q}{q'}\right)\right\rceil+\frac{1}{1-\mu^{\gamma}}+\frac{1}{1-\mu^{{\alpha - \gamma}}}\right).
\end{eqnarray}

\medskip

{\em Case 2}: We now assume that $\|\omega\|_{\infty} = |\omega_{2}| \ge |a_{j}\xi_{1}|$.
Then \eqref{eq:ineq2} applied to \eqref{eq:pf_estimateR2} yields
\begin{eqnarray*}
I_1 &\le& 2^{\gamma''}\sum\limits_{\{j :  |a_{j}\xi_{1}| \le \|\omega\|_{\infty} \}} \min\{|qa_{j}\xi_{1}|^{{\alpha - \gamma}},1\}
\min \{|q'a_{j}\xi_{1}|^{-\gamma},1\} \min \{ \left|q(a_{j}\xi_{1} + 2\omega_{1})\right|^{{\alpha - \gamma}},1\}  \\
& & \qquad \cdot \min\{\left|q'(a_{j}\xi_{1} + 2\omega_{1})\right|^{-\gamma},1\} \sum\limits_{k \in \mathbb{Z}}
\min\left\{\left|\frac rq \left(s_{k}+\frac{\xi_{2}}{a_j^{1/2}\xi_{1}}\right)\right|^{-{\gamma'}},1\right\}\left(\frac rq \right)^{-\gamma''}
\\ &&\qquad \cdot \left|\frac{2\omega_2}{a_j\xi_1} \right|^{-\gamma''}\max\left\{ 1,\left|1+\frac{2\omega_1}{a_j\xi_1} \right|^{\gamma''}\right\}.
\end{eqnarray*}
Applying now \eqref{eq:sk2}, we obtain
\begin{eqnarray}\label{eq:case2} \nonumber
I_1 &\le& 2^{\gamma''}\frac{\left( \frac qr C(\gamma') \right)}{\|2\frac{r}{q}\omega\|_{\infty}^{\gamma''}}\sum_{j \ge 0}
\min\{|qa_{j}\xi_{1}|^{{\alpha - \gamma}},1\}\min
\{|q'a_{j}\xi_{1}|^{-\gamma},1\} \min \{ \left|q(a_{j}\xi_{1} + 2\omega_{1})\right|^{{\alpha - \gamma}},1\}  \\
&& \qquad \cdot \min\{ \left|q'(a_{j}\xi_{1} +
2\omega_{1})\right|^{-\gamma},1\} \left|{a_j\xi_1} \right|^{\gamma''}\max\left\{ 1,\left|1+\frac{2\omega_1}{a_j\xi_1} \right|^{\gamma''}\right\}.
\end{eqnarray}
We next split Case 2 further into the following two subcases.

\medskip

{\em Subcase 2a}: If $1 \le |1+\frac{2\omega_{1}}{a_{j} \xi_{1}}|$, then
\[
\left|{a_j\xi_1} \right|^{\gamma''}\max\left\{ 1,\left|1+\frac{2\omega_1}{a_j\xi_1} \right|^{\gamma''}\right\}
\le \left|a_j\xi_1+2\omega_1 \right|^{\gamma''}.
\]
Hence, by Lemma \ref{lemma:estimate_aj} and exploring inequality \eqref{eq:ineq3},
we can conclude from \eqref{eq:case2} that
\begin{equation}\label{eq:case2a}
I_1 \leq  p\frac{\left( \frac qr C(\gamma') \right)}{\|\frac{q'r}{q}\omega\|_{\infty}^{\gamma''}}
\left( \left\lceil \log_{1/\mu}\left( \frac {q}{q'}\right)\right\rceil+\frac{1}{1-\mu^{\gamma}}+\frac{1}{1-\mu^{{\alpha - \gamma}}}\right).
\end{equation}

\medskip

{\em Subcase 2b}: If $1 \ge |1+\frac{2\omega_{1}}{a_{j}\xi_{1}}|$, then, for all $j \ge 0$,
\[
\min \{ \left|q(a_{j}\xi_{1} + 2\omega_{1})\right|^{{\alpha - \gamma}},1\} \min\{ \left|q'(a_{j}\xi_{1} +
2\omega_{1})\right|^{-\gamma},1\} \max\left\{ 1,\left|1+\frac{2\omega_1}{a_j\xi_1} \right|^{\gamma''}\right\}
\le 1.
\]
Hence, by exploring inequality \eqref{eq:ineq4},
we can conclude from \eqref{eq:case2} that
\begin{equation}\label{eq:case2b}
I_1 \leq  p\frac{\left( \frac qr C(\gamma') \right)}{\|\frac{q'r}{q}\omega\|_{\infty}^{\gamma''}}
\left(\left\lceil \log_{1/\mu}\left( \frac {q}{q'}\right)\right\rceil+\frac{1}{1-\mu^{\gamma'}}+\frac{1}{1-\mu^{{\alpha - \gamma}+\gamma''}}\right).
\end{equation}

We next estimate $I_{2}$. First, notice that the inequality \eqref{eq:pf_estimateR2} still holds for $I_2$ with modified
index set for $j$ -- we have $\{j:|a_j\xi_1| > \|\omega\|_{\infty} \}$ instead of $\{j:|a_j\xi_1| \le \|\omega\|_{\infty} \}$ in this case.
Therefore, by \eqref{eq:sk2},
\begin{eqnarray*}
I_2 &\leq&
\sum\limits_{\{j :  |a_{j}\xi_{1}| > \|\omega\|_{\infty} \}} \min\{|qa_{j}\xi_{1}|^{{\alpha - \gamma}},1\}\min
\{|q'a_{j}\xi_{1}|^{-\gamma},1\} \sum\limits_{k \in \mathbb{Z}}
\min\left\{\left|\frac rq \left(s_{k}+\frac{\xi_{2}}{a_j^{1/2}\xi_{1}}\right)\right|^{-{\gamma}}\hspace*{-0.25cm},1\right\} \\
&\le&
\left( \frac qr C(\gamma) \right)\sum_{\{j:|a_j\xi_1| > \|\omega\|_{\infty} \}} \min\{|qa_{j}\xi_{1}|^{{\alpha - \gamma}},1\}
\min \{|q'a_{j}\xi_{1}|^{-\gamma},1\} \\
&\le& \left( \frac qr C(\gamma)\right) \sum_{\{j:|a_j\xi_1| > \|\omega\|_{\infty} \}} |q'a_j\xi_1|^{-\gamma}
\end{eqnarray*}
By Lemma \ref{lemma:estimate_aj}, this estimation can be finalized to
\begin{equation}\label{eq:estimatei2}
I_2 \le p\frac{\frac qr C(\gamma)}{\|q'\omega\|_{\infty}^{\gamma}}\Bigl( \frac{1}{1-\mu^{\gamma}} \Bigr).
\end{equation}

\medskip

We are now ready to prove the claimed estimate for $R(c)$ using \eqref{eq:pf_estimateR1}. For this, we define the constants $T_1, T_2$, and $T_3$ as in
Proposition \ref{prop:estimateR}. Now define
\[
{\mathcal Q} = \{m \in \ZZ^2 : |m_1| > |m_2|\} \quad \text{and} \quad \tilde{{\mathcal Q}} = \{m \in \ZZ^2 : c_1^{-1}|m_1|
 > c_2^{-1}|m_2|\}.
\]
If $c_1^{-1}|m_1| > c_2^{-1}|m_2|$, then, by \eqref{eq:case1} and \eqref{eq:estimatei2},
\[
\Gamma(\pm M_c^{-1}m) \le T_1\|m\|^{-\gamma}_{\infty}+T_3\|m\|^{-\gamma}_{\infty} \quad \text{for all}\quad m \in \tilde{{\mathcal Q}}.
\]
If $c_1^{-1}|m_1| \le c_2^{-1}|m_2|$ with $m \neq 0$, then, by \eqref{eq:case2a}, \eqref{eq:case2b}, and \eqref{eq:estimatei2},
\[
\Gamma(\pm M_c^{-1}m) \le T_2\|m\|^{-\gamma''}_{\infty}+T_3\|m\|^{-\gamma}_{\infty} \quad
\text{for all}\quad m \in \tilde{{\mathcal Q}}^c\backslash \{0\}.
\]
Therefore, we obtain
\begin{eqnarray}\nonumber
R(c) &=& \sum_{m \in \ZZ^2\backslash\{0\}}\left( \Gamma(M_c^{-1}m)\Gamma(-M_c^{-1}m)\right)^{1/2}\\ \label{eq:Rc}
&\le& \left(\sum_{m \in \tilde{{\mathcal Q}}} T_1\|m\|^{-\gamma}_{\infty}+T_3\|m\|^{-\gamma}_{\infty}\right) +
\left(\sum_{m \in \tilde{{\mathcal Q}}^c\backslash\{0\}} T_2\|m\|^{-\gamma''}_{\infty}+T_3\|m\|^{-\gamma}_{\infty}\right)
\end{eqnarray}
Notice that, since $\tilde{{\mathcal Q}} \subset {\mathcal Q}$,
\[
\sum_{m \in \tilde{{\mathcal Q}}} \|m\|^{-\gamma} \le \sum_{m \in {{\mathcal Q}}} \|m\|^{-\gamma}.
\]
Also, we have
\[
\sum_{m \in \tilde{{\mathcal Q}}^c\backslash\{0\}} \|m\|^{-\gamma''} \le
\min \left \{\left \lceil \frac{c_1}{c_2} \right \rceil,2 \right \} \sum_{m \in {{\mathcal Q}^c\backslash\{0\}}} \|m\|^{-\gamma''}.
\]
Therefore, \eqref{eq:Rc} can be continued by
\begin{equation*}
R(c) \le T_3 \sum_{m \in \ZZ^2\backslash\{0\}}\|m\|_{\infty}^{-\gamma} + T_1 \sum_{m \in {\mathcal Q}} \|m\|_{\infty}^{-\gamma}
+ \min \left \{ \left \lceil \frac{c_1}{c_2} \right \rceil, 2 \right \} T_2 \sum_{m \in {{\mathcal Q}^c\backslash\{0\}}} \|m\|^{-\gamma''}.
\end{equation*}
}
To provide an explicit estimate for the upper bound of $R(c)$, we
compute $\sum_{m \in {\mathcal Q}} \|m\|_{\infty}^{-\gamma}$ and $\sum_{m \in {\mathcal Q}^c\backslash} \|m\|_{\infty}^{-\gamma}$ as follows:
\begin{eqnarray*}
\sum_{m \in {\mathcal Q}} \|m\|^{-\gamma}_{\infty} &=& 2\sum_{m_1=1}^{\infty}|m_1|^{-\gamma}+4\sum_{m_2=1}^{\infty}\sum_{m_1=m_2+1}^{\infty}|m_1|^{-\gamma}
\\
&\leq& 2\Bigl(1+\frac{1}{\gamma-1}\Bigl)+\frac{4}{\gamma-1}\Bigl(1+\frac{1}{\gamma-2}\Bigr).
\end{eqnarray*}
and
\begin{eqnarray*}
\sum_{m \in {\mathcal Q}^c\backslash \{0\}} \|m\|^{-\gamma}_{\infty}
&=&\sum_{m \in {\mathcal Q}} \|m\|^{-\gamma}_{\infty}+4\sum_{m_1=1}^{\infty}|m_1|^{-\gamma}\\
&\leq& 6\Bigl(1+\frac{1}{\gamma-1}\Bigl)+\frac{4}{\gamma-1}\Bigl(1+\frac{1}{\gamma-2}\Bigr)
\end{eqnarray*}
This completes the proof.

\subsection{Proofs of Results from Section \ref{subsec:compact}}
\label{subsec:compactproof}
\subsubsection{Proof of Lemma \ref{lemma:m0}}
\label{subsec:proofm0}
First, it is obvious that the function $|m_0|^2$ is even. Next,
letting $y = \sin^2(\pi\xi_1)$, the values of the trigonometric polynomial $|m_0(\xi_1)|^2$ can be expressed in the form
\[
|m_0(\xi_1)|^2 = P(y), \quad \mbox{where } P(y) = (1-y)^K\sum_{n=0}^{L-1}{K-1+n\choose n}y^n.
\]
We compute
\begin{eqnarray*}
\lefteqn{P'(y)}\\
&=& \hspace*{-0.1cm} (1-y)^{K-1}\Bigl[-K\sum_{n=0}^{L-1}{K-1+n \choose n}y^n+(1-y)\sum_{n=1}^{L-1}n{K-1+n \choose n}y^{n-1}\Bigr] \\
&=& \hspace*{-0.1cm} (1-y)^{K-1}\Bigl[-K\sum_{n=0}^{L-1}{K-1+n \choose n}y^n-\sum_{n=1}^{L-1}n{K+n-1\choose n}y^n
+\sum_{n=0}^{L-2}(n+1){K+n\choose n+1}y^n\Bigr] \\
&=& \hspace*{-0.1cm} (1-y)^{K-1}\Bigl[-K{K+L-2\choose L-1}y^{L-1}-(L-1){K+L-2\choose L-1}y^{L-1}\Bigr] \\
&=& \hspace*{-0.1cm} -(K+L-1){K+L-2\choose L-1}y^{L-1}(1-y)^{K-1}.
\end{eqnarray*}
Hence, $P'(y) < 0$ for $y \in (0,1)$, and this immediately implies that $|m_0(\xi_1)|^2$ is decreasing on $(0,\frac12)$.
The second derivative can be derived by
\begin{eqnarray*}
\lefteqn{P''(y)}\\
&=& -(L-1)(K+L-1){K+L-2\choose L-1}y^{L-2}(1-y)^{K-1}+(K+L-1)(K-1) \\
&& \qquad \cdot {K+L-2\choose L-1}y^{L-1}(1-y)^{K-2} \\
&=& (K-1)(K+L-1)\Bigl(\frac{K+L-2}{K-1}\Bigr){K+L-2\choose L-1}y^{L-2}(1-y)^{K-2}\Bigl( y-\frac{L-1}{K+L-2}\Bigr).
\end{eqnarray*}
Thus $P''(y) < 0$ for $0 < y < \frac{L-1}{K+L-2}$. Hence
\[
\frac{\partial^2}{\partial \xi_1^2}|m_0(\xi_1)|^2 < 0 \quad \mbox{for }
0<\xi_1<\frac{1}{\pi}\arcsin\Bigl[ \Bigl(\frac{L-1}{K+L-2} \Bigr)^{\frac 12}\Bigr].
\]
Since $\frac{L-1}{K+L-2} \ge \frac{1}{4}$ implies
\[
\frac{1}{\pi}\arcsin\Bigl[ \Bigl(\frac{L-1}{K+L-2} \Bigr)^{\frac 12}\Bigr] \ge \frac 16,
\]
concavity of $|m_0|^2$ on $(0,\frac16)$ is proven.

\subsubsection{Proof of Lemma \ref{lemma:tildem}}
\label{subsec:prooftildem}

Letting $y = \sin^2(\pi\xi_1)$ for $\xi_1 \in [0,1]$, the values of the trigonometric polynomial $|\tilde{m}_0(\xi_1)|^2$
can be expressed in the form
\[
|\tilde{m}_0(\xi_1)|^2 = \tilde{P}(y), \quad \mbox{where } \tilde{P}(y) = (1-y)^{K'}\sum_{n=0}^{L-1}{K-1+n\choose n}y^n.
\]
If $K'=0$, $\tilde{P}(y)$ is increasing on $(0,1)$, which proves one part of claim (i).

To prove the remaining part of claim (i), let us assume that $K'>0$. For $y \in (0,1)$, direct computation yields
\begin{eqnarray*}
\lefteqn{(\tilde{P})'(y)}\\ &=& (1-y)^{K'-1}\Bigl[(K-K')\sum_{n=0}^{L-2}{K-1+n\choose n}y^{n}-(K'+L-1){K-L-2\choose L-1}y^{L-1}\Bigr] \\
&=& y^{L-1}(1-y)^{K'-1}\Bigl[(K-K')\sum_{n=0}^{L-2}{K-1+n\choose n}y^{n-L+1}-(K'+L-1){K-L-2\choose L-1}\Bigr].
\end{eqnarray*}
Hence, $(\tilde{P})'(y)$ is a product of a decreasing function and of $y^{L-1}(1-y)^{K'-1}>0$ on $(0,1)$.
This implies that $\tilde{P}(y)$ is increasing on $(0,z_0)$ as long as $(\tilde{P})'(z_0) \ge 0$.
Obviously, the same is true for $|\tilde{m}_0(\xi_1)|^2$. Thus, since $\frac 14 = \sin^2(\frac{\pi}{6})$,
it suffices to show that $(\tilde{P})'(\frac 14) \ge 0$.  For this, we compute
\[
\tilde{P}'(\tfrac14) = (K'+L-1)\Bigl(\frac{3}{4}\Bigr)^{K'-1}4^{1-L}\Bigl[\frac{K-K'}{K'+L-1}\sum_{n=0}^{L-2}{K-1+n\choose n}2^{2L-2}2^{-2n}-{K+L-2\choose L-1}\Bigr].
\]
Since $K \ge 7$ and $\frac{K+L-2}{L-1} \le 4$, the term in the bracket $[ \:\cdot\: ]$ can be further estimated by
\begin{eqnarray*}
2^{2L}\Bigl(\frac{K-K'}{K'+L-1}\Bigr)\Bigl(\frac{1}{4}\Bigr)\sum_{n=0}^{L-2}{K-1+n\choose n}2^{-2n}-{K+L-2\choose L-1}
&&\\
&& \hspace{-300pt}  \ge 2^{2L} \Bigl(\frac{K-K'}{K'+L-1}\Bigr)\Bigl(\frac{1}{4}\Bigr)\Bigl[ 1+\frac{K}{4}+\frac{K(K+1)}{4^2 \cdot 2}+
\frac{K(K+1)(K+2)}{4^3 \cdot 6}\Bigr]-{K+L-2\choose L-1} \\
&& \hspace{-300pt} \ge 2^{2L}\Bigl(\frac{K-K'}{K'+L-1}\Bigr)- {K+L-2 \choose L-1}\\
&& \hspace{-300pt} = 2^{2L}\Bigl( \frac{K-K'}{K'+L-1}\Bigr)-\frac{(K+L-2)\cdot \ldots \cdot K)}{(L-1)!}
\ge 2^{2L}\Bigl( \frac{K-K'}{K'+L-1}\Bigr)-2^{2L-2} \ge 0,
\end{eqnarray*}
where the last inequality follows from $\frac{K-K'}{K'+L-1} \ge \frac 14$. This completes the proof of claim (i).

To prove claim (ii), observe that, for $0\le y\le1$,
\begin{eqnarray*}
|\tilde{P}(y)| &\leq& \sum_{n=0}^{L-1}{K-1+n\choose n} \max_{y \in [0,1]} |(1-y)^{K'}||y^n|\\
&=& \sum_{n=0}^{L-1}{K-1+n\choose n}\Bigl( \frac{K'}{K'+n}\Bigr)^{K'}\Bigl(\frac{n}{K'+n}\Bigr)^{n}.
\end{eqnarray*}
Thus claim (ii) is proven.

\subsubsection{Proof of Lemma \ref{lemma:gamma}}

Without loss of generality we assume that $L$ is even, since the `odd case' can be proven
similarly. Since $K'=0$ and by obvious rules for binomial coefficients,
\[
C_2 = \sum_{n=0}^{L-1}{K-1+n\choose n} = \sum_{n=0}^{K-1}{K-1+n\choose n}-\sum_{n=L}^{K-1}{K-1+n\choose n}
\]
Utilizing the estimate $\sum_{n=0}^{K-1}{K-1+n\choose n} < 2^{2K-2}$ from \cite{Dau92}, we have
\begin{equation}\label{eq:eq2}
C_2 < 2^{2K-2}-\sum_{n=L}^{K-1}{K-1+n\choose n}
\end{equation}
Now let $a_n = {K-1+n\choose n}$ for $n = 0,\dots,L-1$, and note that
\[
\frac{a_{n+1}}{a_n} = \frac{K+n}{n+1} \ge \frac{K+L-1}{L} \ge 2,
\]
which implies
\begin{eqnarray*}
a_L+\dots+a_{K-1} &\ge& a_L+a_{L+1}+\dots+a_{\frac{3L}{2}-1} \\ &\ge& 2^{L-2}(2a_1)+2^{L-3}(2a_3)+\dots+2^{L/2-1}(2a_{L-1})\\
&\ge& 2^{L/2-1}(a_0+\dots+a_{L-1}).
\end{eqnarray*}
Thus, using \eqref{eq:eq2},
\[
C_2 <2^{2K-2}-2^{L/2-1}C_2,
\]
and hence
\[
C_2 < 2^{2K-2-L/2+1} = 2^{2K-L/2-1}.
\]
Thus, we have
\[
\max_{\xi_1 \in \bR} |\tilde{m}_0(\xi_1)|^2 \le 2^{2K-L/2-1}
\]
and
\[
2\gamma = 2K-\log_2\left(C_2\right) > \frac{L}{2}+1,
\]
which is what was claimed.

\subsubsection{Proof of Proposition \ref{proposition:upperlowershear}}
By (i) and (ii) in Lemma \ref{lemma:m0}, we obtain
\begin{equation}\label{eq:lowerm1}
|m_1(4\xi_1)|^2 \ge |m_0(\tfrac16)|^2\chi_{[1/12,1/6]\cup[-1/12,-1/6]}(\xi_1).
\end{equation}
On the other hand, by definition of $m_1$,
\begin{eqnarray*}
|m_1(\xi_1)|^2 &=& (\sin(\pi\xi_1))^{2K}\sum_{n=0}^{L-1}{K-1+n\choose n}(\cos(\pi\xi_1))^{2n} \nonumber \\
&=& (\sin(\pi\xi_1))^{2(K-K')}(\sin(\pi\xi_1))^{2K'}\sum_{n=0}^{L-1}{K-1+n\choose n}(\cos(\pi\xi_1))^{2n} \nonumber \\
&\le& |\pi\xi_1|^{2(K-K')}\sum_{n=0}^{L-1}{K-1+n\choose n}\left( \frac{K'}{K'+n}\right)^{K'}\left( \frac{n}{K'+n}\right)^n \nonumber \\
&\le& C_2|\pi\xi_1|^{2(K-K')}.
\end{eqnarray*}
Therefore,
\begin{equation}\label{eq:upperm1}
|m_1(\xi_1)|^2 \le \min\{1,C_2|\pi\xi_1|^{2(K-K')}\}
\end{equation}
since $|m_1(\xi_1)|^2 \le 1.$
Let now
\[
q = 4\pi (C_2)^{1/(2(K-K'))}, \quad q' = 2\pi \left(C_2\prod_{j=0}^{J_1-1}\left|\tilde{m}_0\left(\frac{2^{-j}}{2\pi}\right)\right|^2
e^{2^{-J_1+1}\sum_{n} |h(n)||n|}\right)^{-\frac{1}{2\gamma}}, \quad \mbox{and } r=2q',
\]
with $C_2$ being defined in \eqref{eq:c2}. It is easy to check that $q > r > q' > 0$.
By Propositions \ref{proposition:piLower} and \ref{proposition:piUpper},\eqref{eq:lowerm1} and \eqref{eq:upperm1},
we obtain \eqref{eq:uppershearlet} and \eqref{eq:lowershearlet}.

In \eqref{eq:uppershearlet}, the decay rate $\gamma$ is given by $\gamma = K-K'-\frac 12\log_2(C_2)$. Hence, by Lemma \ref{lemma:gamma},
in the special case that $K'=0$ and $L \ge 10$, we have $\gamma > L/4+1/2 \ge 3$. This implies the feasibility
condition \eqref{eq:psi}.


\end{document}